\documentclass{article}

\usepackage{proof}
\usepackage{latexsym}
\usepackage{amssymb}
\usepackage{cite}

\newcommand{\blem}{\begin{lemma}}
\newcommand{\elem}{\end{lemma}}
\newcommand{\bth}{\begin{theorem}}
\newcommand{\ethm}{\end{theorem}}
\newcommand{\benu}{\begin{enumerate}}
\newcommand{\eenu}{\end{enumerate}}
\newcommand{\bdes}{\begin{description}}
\newcommand{\edes}{\end{description}}
\newcommand{\bdf}{\begin{definition}}
\newcommand{\edf}{\end{definition}}
\newcommand{\bcor}{\begin{cor}}
\newcommand{\ecor}{\end{cor}}
\newcommand{\bprp}{\begin{proposition}}
\newcommand{\eprp}{\end{proposition}}
\newcommand{\bmlem}{\begin{mlemma}}
\newcommand{\emlem}{\end{mlemma}}
\newcommand{\bclm}{\begin{claim}}
\newcommand{\eclm}{\end{claim}}
\newcommand{\bprf}{{\bf Proof}.\hspace{2mm}}

\newcommand{\eprf}{\hspace*{\fill} $\Box$}

\newcommand{\beqn}{\begin{equation}}
\newcommand{\eeqn}{\end{equation}}
\newcommand{\beqnarr}{\begin{eqnarray}}
\newcommand{\eeqnarr}{\end{eqnarray}}
\newcommand{\beqnarrs}{\begin{eqnarray*}}
\newcommand{\eeqnarrs}{\end{eqnarray*}}

\newcommand{\spand}{\,\&\,}

\newcommand{\restrict}{\!\upharpoonright\!}

\newtheorem{theorem}{Theorem}[section]
\newtheorem{definition}[theorem]{Definition}
\newtheorem{proposition}[theorem]{Proposition}
\newtheorem{lemma}[theorem]{Lemma}
\newtheorem{cor}[theorem]{Corollary}
\newtheorem{mlemma}[theorem]{Main Lemma}
\newtheorem{claim}[theorem]{Claim}

\newcommand{\alp}{\alpha}

\newcommand{\del}{\delta}
\newcommand{\Del}{\Delta}
\newcommand{\ome}{\omega}

\newcommand{\bet}{\beta}
\newcommand{\gam}{\gamma}
\newcommand{\Gam}{\Gamma}
\newcommand{\kap}{\kappa}
\newcommand{\sig}{\sigma}

\newcommand{\Lam}{\Lambda}
\newcommand{\vphi}{\varphi}

\newcommand{\fal}{\forall}
\newcommand{\exi}{\exists}

\newcommand{\Rarw }{\Rightarrow}

\newcommand{\lrarw}{\leftrightarrow}
\newcommand{\Lrarw}{\Leftrightarrow}

\newcommand{\calL}{{\cal L}}

\newcommand{\calP}{{\cal P}}

\newcommand{\st}{\star}

\newcommand{\la}{\langle}
\newcommand{\ra}{\rangle}

\newcommand{\msfiv}{\mbox{\hspace{5mm}}}

\title{Proof search in multi-succedent sequent calculi for intuitionistic logic under the Sch\"utte's schema
}
\author{Toshiyasu Arai
\\
Graduate School of Science,
Chiba University
\\
1-33, Yayoi-cho, Inage-ku,
Chiba, 263-8522, JAPAN
\\
tosarai@faculty.chiba-u.jp
}
\date{}
\begin{document}
\maketitle

\begin{abstract}
In this note
we present terminating and bicomplete proof searches in multi-succedent sequent calculi for intuitionistic propositional 
logic,
fragments of intuitionistic predicate logic
and full intuitionistic predicate logic
in the spirit of Sch\"utte's schema.
\end{abstract}

In this note we present terminating and bicomplete proof search procedures
in multi-succedent sequent calculi for intuitionistic propositional logic and fragments of intuitionistic predicate logic 
with respect to the Kripke semantics.
G. Mints\cite{Mintsint, Mintsdraft} in his later years investigated a proof search procedure in single-succedent sequent calculus for 
intuitionistic predicate logic.

\section{Sch\"utte's schema}\label{sec:prfsearchLJm}

\textit{Sch\"utte's schema} in \cite{Schuette56} is described as follows.

Given a logic calculus, e.g., a sequent calculus and a semantics for the logic,
search recursively a cut-free derivation of a given sequent in a bottom-up manner.
This results in a (computable but possibly infinite) deduction tree of the given sequent.
If the tree is a (finite) derivation, then it tells us that the sequent is cut-free derivable in the sequent calculus.
Otherwise it yields a counter model of the sequent with respect to the semantics.
Let us call the division \textit{Sch\"utte's dichotomy}.
Thus the schema shows simultaneously the completeness of the (cut-free fragment of) the sequent calculus
with respect to the semantics and the Hauptsatz for the calculus.
The schema has been successfully applied to (first-order and higher-order) classical logic calculi
by K. Sch\"utte\cite{Schuette60}.

G. Mints\cite{Mintsint, Mintsdraft} investigates a proof search procedure in single-succedent sequent calculus for 
intuitionistic predicate logic with respect to Kripke semantics\cite{Kripke}.
On the other side a multi-succedent sequent calculus for intuitionistic logic was introduced in \cite{Maehara},
and it is known that it relates to semantic tableaux and Kripke semantics, cf.\,\cite{TroSch}.
However a naive proof search procedure in a multi-succedent sequent calculus for intuitionistic propositional logic
may not terminate in a finite number of steps since the left rule for implication might be iterated unlimitedly.
For example the following proof search for the sequent sequent $(p\supset q)\supset \bot\Rarw q$ does not terminate.
\[
\infer{(p\supset q)\supset \bot\Rarw q}
{
 \infer{(p\supset q)\supset \bot\Rarw p\supset q}
 {
 \infer{(p\supset q)\supset \bot,p\Rarw q}
  {
  \infer{(p\supset q)\supset \bot,p\Rarw p\supset q}
    {\infer*{(p\supset q)\supset \bot,p\Rarw q}{}}
  &
   \bot,p \Rarw q
   }
  }
&
\bot\Rarw q
}
\]
It is desirable that a proof search procedure for a decidable logic is terminating.

In the second part of this note let us give terminating and bicomplete proof search procedures under the schema
in multi-succedent sequent calculi ${\sf LJpm}^{{\tt c}}$ and ${\sf LJm}^{{\tt c}}$ 
for intuitionistic propositional logic {\sf Ip} and fragments of intuitionistic predicate logic {\sf Iq}
with respect to the Kripke semantics.
By a bicomplete proof search procedure we mean a procedure such that we can extract a counter model
from failed proof search, cf.\,\cite{Dyckhoff}.
${\sf LJpm}^{{\tt c}}$ denotes a modified multi-succedent sequent calculus for intuitionistic propositional logic.
In the calculus conclusions of minor formulas are added in succedents when a right rule is applied.
It turns out that a proof search procedure for ${\sf LJpm}^{{\tt c}}$ is terminating.
A terminating proof search for the above sequent $(p\supset q)\supset \bot\Rarw q$ in ${\sf LJpm}^{{\tt c}}$ runs as follows.
\[
\infer[(L\supset)]{(p\supset q)\supset \bot\Rarw q}
{
 \infer{(p\supset q)\supset \bot\Rarw p\supset q,q}
 {
 \infer[(L\supset)]{(p\supset q)\supset \bot,p\Rarw q}
  {
  (p\supset q)\supset \bot,p\Rarw p\supset q,q
  &
   \bot,p \Rarw q
   }
  }
&
\bot\Rarw q
}
\]
When the left rules $(L\supset)$ are applied with the main formula $(p\supset q)\supset \bot$,
the left minor formula $p\supset q$ is accompanied with its conclusion $q$.
The leaf sequent $(p\supset q)\supset \bot,p\Rarw p\supset q,q$ is not a classical tautology.

In Section \ref{sec:Ip} we consider a proof search in a sequent calculus ${\sf LJpm}^{{\tt c}}$
for intuitionistic propositional logic {\sf Ip}.
A related work is done in \cite{PintoDyckhoff}.

Given a sequent $S_{0}$, our search procedure yields a finite $(\land,\lor)$-tree $TR(S_{0})$ 
of finite deductions.
Each leaf in the $(\land,\lor)$-tree receives a value in $\{0,1\}$.
Any formula in the sequent at a topmost node is essentially either an atom (propositional variable) or an absurdity $\bot$, 
and
the antecednt contains no $\bot$ and there is no common atom in the antecedent and the succedent when the value is $0$.
Otherwise the sequent is an axiom.
Furthermore if the value of the whole $(\land,\lor)$-tree is $1$, then by pruning,
we can extract a derivation of $S_{0}$
from the tree of deductions $TR(S_{0})$.
Otherwise a Kripke model is readily constructed, in which $S_{0}$ is falsified.

As a corollary and by the depth-first left-to-right implementation of proof search,
we see in Section \ref{subsec:IpPSPACE}
that the intuitionistic propositional logic is in PSPACE.
The fact was first proved by R. Ladner\cite{Ladner}.

In Section \ref{sec:Iq}  we consider proof search procedures in sequent calculi ${\sf LJm}^{{\tt c}}$
for fragments of the intuitionistic predicate logic {\sf Iq}.
Proof search procedures for the fragments are terminating and bicomplete.
As a corollary we see that the fragments are decidable.
One fragment denoted by $\fal^{(++)}\exi^{(-)}$ is a subclass of the positive fragment of intuitionistic predicate logic.
Mints\cite{MintsSol} showed that the positive fragment is decidable.
It turns out that the fragment $\fal^{(++)}\exi^{(-)}$ is is in PSPACE.
The other fragment $\fal^{(++,-)}\exi^{(+)}$ is shown to be solvable in exponential space.

In Section \ref{sec:Iqfull}
we consider a proof search procedure in the sequent calculi  {\sf LJm} and ${\sf LJm}+(cut)$ with the cut rule.
for the full intuitionistic predicate logic {\sf Iq}.
Proof search procedures yields a (possibly infinite) tree of finite deduction trees for a given sequent.
Here the search tree may be infinite since we need to test all of terms in the infinite list.
In defining Kripke models from the tree of finite deductions, we may encounter an inconsistency.
To avoid the inconsistency in the definition of Kripke models, we need to transform the tree of deductions.
However if we need only to show the completeness of {\sf LJm} with the cut rule $(cut)$, 
then it turns out in subsection \ref{subsec:completeness}
that the inconsistency disappears, and we don't need to transform trees of deductions.
This is done by using characteristic formulas as in \cite{Kripke}.
In subsection \ref{subsec:transfer}
we consider a transformation of trees of deductions inspired by the transfer rule in \cite{Kripke}
to avoid an inconsistency in the definition of Kripke models.
From the transformed tree of deductions, we obtain a Sch\"utte's dichotomy, cf.\,Theorem \ref{th:schuttedichotomy}.

\section{Propositional case}\label{sec:Ip}

In this section we consider a proof search procedure in 
a sequent calculus ${\sf LJpm}^{{\tt c}}$ for the intuitionistic propositional logic {\sf Ip}.
In the calculus ${\sf LJpm}^{{\tt c}}$ the derived objects are 
\textit{sequents}, which are ordered pairs of cedents denoted by $\Gam\Rarw\Del$,
where $\Gam$ is the \textit{antecedent} and $\Del$ the \textit{succedent} of the sequent.
Though the calculus is for intuitionistic propositional logic {\sf Ip},
several formulas may occur in succedents, i.e., multi-succedent sequents.
{\sf p} in ${\sf LJpm}^{{\tt c}}$ stands for \textit{propositional},
{\sf m} for \textit{multi-succedent sequents}.
Moreover the super script {\tt c} in ${\sf LJpm}^{{\tt c}}$ indicates \textit{conclusion}.
The conclusion $\bet$ of an implicational minor (active) formula $\alp\supset\bet$ in succedents
is augmented in each inference rule.
For example a right rule for disjunction in ${\sf LJpm}^{{\tt c}}$ is of the form
\[
\infer{\Gam\Rarw\Del,(p\supset q)\lor r}
{
\Gam\Rarw \Del,p\supset q, q
}
\]
Let us explain the reason why we augment conclusions in succedents.
Let {\sf LJpm} be a sequent calculus for the intuitionistic propositional logic {\sf Ip}
in which conclusions may be absent in succedents.
Namely {\sf LJpm} is the propositional fragment of the calculus {\sf m-G3i} in \cite{TroSch}.
A right rule for disjunction in {\sf LJpm} is of the form
\[
\infer{\Gam\Rarw\Del,(p\supset q)\lor r}
{
\Gam\Rarw \Del,p\supset q
}
\]
In searching a derivation of a given sequent,
we need to analyze implicational formulas $\alp\supset\bet$ in antecedents several times.
\[
\infer[(L\supset)]{\alp\supset\bet,\Gam\Rarw\Del}
{
\alp\supset\bet,\Gam\Rarw\Del,\alp
&
\bet,\Gam\Rarw\Del
}
\]
Although we can control the number of applications of the left rule $(L\supset)$
in terms of the given sequents as in \cite{BussIemhoff},
let us take another route.
By augmenting conclusions in succedents, 
we see that each formula $\alp\supset\bet$ is used as the major formula of a right rule $(R\supset)$ at most once
on each branch in the searching tree.
Moreover each formula $\alp\supset\bet$ is used as the major formula of a left rule $(L\supset)$ at most once
in a part of a branch such that the part contains no right rule $(R\supset)$, cf.\,Corollary \ref{prp:depthIp}.
In this way we can conclude that the search procedure terminates, cf.\, Lemma \ref{prp:termination}.
\\

A formula of the form $\alp\supset\bet$ is an \textit{implicational} formula.

Let $\bet_{i}$ be formulas.
$\bet_{1}\supset\bet_{2}\supset\cdots\supset\bet_{n}\supset\bet_{n+1}$
denotes the formula $\bet_{1}\supset(\bet_{2}\supset(\cdots\supset(\bet_{n}\supset\bet_{n+1})\cdots))$
in the association to the right.

\bdf\label{df:c}
{\rm 

Let $\alp\equiv(\bet_{1}\supset\bet_{2}\supset\cdots\supset\bet_{n}\supset\bet_{n+1})\,(n\geq 0)$
with a non-implicational formula $\bet_{n+1}$, and
$\alp_{k}\equiv(\bet_{k}\supset\bet_{k+1}\supset\cdots\supset\bet_{n}\supset\bet_{n+1})$ for $k=1,2,\ldots,n$.
Then let
$\alp^{{\tt c}}=\{\bet_{n+1}\}$
and 
$\alp^{{\tt p}}=\{\bet_{k}: k=1,2,\ldots,n\}$.
}
\edf

Note that 
$\alp\lor\alp^{{\tt c}}$ is intuitionistically equivalent to $\alp$.

\subsection{Sequent calculus ${\sf LJpm}^{{\tt c}}$ for the intuitionistic propositional logic {\sf Ip}}\label{subsec:LJpm}
{\bf Axioms}.
\[
(T)\: \Gam\Rarw \Del  \mbox{ if } \Gam\cap\Del\cap Atm\neq\emptyset \msfiv
(\bot)\: \Gam\Rarw\Del \mbox{ if } \bot\in\Gam
\]
where $Atm$ denotes the set of atoms.
{\bf Inference rules}.
\[
\infer[(L\lor)]{\Gam\Rarw\Del}
{\alp_{0},\Gam\Rarw\Del
&
\alp_{1},\Gam\Rarw\Del
}
\mbox{ with } (\alp_{0}\lor\alp_{1})\in\Gam
\]
\[
\infer[(R\lor)]{\Gam\Rarw\Del}
{
\Gam\Rarw\Del,\alp_{0},\alp_{0}^{{\tt c}},\alp_{1},\alp_{1}^{{\tt c}}
}
\mbox{ with } (\alp_{0}\lor\alp_{1})\in\Del
\]
\[
\infer[(L\land)]{\Gam\Rarw\Del}
{
\alp_{0},\alp_{1},\Gam\Rarw\Del
}
\mbox{ with } (\alp_{0}\land\alp_{1})\in\Gam
\]
\[
\infer[(R\land)]{\Gam\Rarw\Del}
{
\Gam\Rarw\Del,\alp_{0},\alp_{0}^{{\tt c}}
&
\Gam\Rarw\Del,\alp_{1},\alp_{1}^{{\tt c}}
}
\mbox{ with } (\alp_{0}\land\alp_{1})\in\Del
\]
\[
\infer[(L\supset)]{\Gam\Rarw\Del}
{
\{
\Gam\Rarw\Del,\bet,\bet^{{\tt c}}
: \bet\in\alp^{{\tt p}}\}
&
\alp^{{\tt c}},\Gam\Rarw\Del
}
\mbox{ with an implicational } \alp\in\Gam
\]
\[
\infer[(R\supset)]{\Gam\Rarw\Del}
{\alp^{{\tt p}},\Gam\Rarw\alp^{{\tt c}}}
\mbox{ with an implicational } \alp\in\Del
\]

Each sequent above the line is an \textit{upper sequent}, and the sequent below the line is the \textit{lower sequent} of an inference rule.
For example in a $(L\lor)$, $\alp_{i},\Gam\Rarw\Del$ is an upper sequent for $i\in\{0,1\}$, and $\Gam\Rarw\Del$ is the lower sequent.
Observe that the antecedent of the lower sequent is a subset of the antecedent of each upper sequent
in any inference rules, and the succedent of the lower sequent is a subset of the succcedent of each upper sequent
in an inference rule other than $(R\supset)$.

\textit{Derivations} in ${\sf LJpm}^{{\tt c}}$ are defined as usual.
These are labelled finite trees whose leaves are labelled axioms and which are locally correct with respect to inference rules.
While \textit{deductions} are labelled finite trees which are locally correct with respect to inference rules.
In deductions labels of leaves may be any sequents.
\\

A \textit{Kripke frame} is a quasi order $\langle W, \preceq\rangle$.
This means that $W\neq\emptyset$ and $\preceq$ is a reflexive and transitive relation on $W$.
A \textit{Kripke model} is a triple $\langle W,\preceq,V\rangle$, where $\langle W, \preceq\rangle$ is a Kripke frame,
and $V:W\to \calP(Atm)$ such that $V(w)\subset V(v)$ if $w\preceq v$.

For formulas $\alp$ and $w\in W$, $w\models\alp$ is defined recursively.

\benu
\item
$w\models p$ iff $p\in V(w)$.
\item
$w\models\alp\lor\bet$ iff $w\models\alp$ or $w\models\bet$.
\item
$w\models\alp\land\bet$ iff $w\models\alp$ and $w\models\bet$.
\item
$w\not\models\bot$.
\item
$w\models(\alp\supset\bet)$ iff for every $v\succeq w$, if $v\models\alp$, then $v\models\bet$.
\eenu

A sequent $\Gam\Rarw\Del$ is \textit{intuitionistically valid} if 
$w\models \bigwedge\Gam$ implies $w\models\bigvee\Del$
for any Kripke model $\langle W,\preceq,V\rangle$ and any $w\in W$.

\bprp
Any derivable sequent in ${\sf LJpm}^{{\tt c}}$ is intuitionistically valid.
\eprp

\subsection{Proof search in ${\sf LJpm}^{{\tt c}}$}\label{subsec:prfsearchLJpm}

It is easy to define a terminating and bicomplete
proof search procedure in a sequent calculus {\sf LKp} for the classical propositional logic, cf.\,\cite{Mintsmodal}.
Actually the procedure yields a conjunctive normal form 
$\bigwedge_{i}\{p_{i1},\ldots,p_{in_{i}}\Rarw q_{i1},\ldots,q_{im_{i}}\}$ of a given formula,
where each $p_{i1},\ldots,p_{in_{i}}\Rarw q_{i1},\ldots,q_{im_{i}}$ is obtained from a leaf sequent in the search tree
by deleting all of non-atomic formulas.
The sequent $p_{i1},\ldots,p_{in_{i}}\Rarw q_{i1},\ldots,q_{im_{i}}$ is classically equivalent to
the disjunctive formula $\lnot p_{i1}\lor\cdots\lor\lnot p_{in_{i}}\lor q_{i1}\lor\cdots\lor q_{im_{i}}$.

On the other side each $w\in W$ in a Kripke model $\la W,\preceq,V\ra$ determines a classical truth assignment
which assigns a truth value to each atom and to each implicational formula.
The truth value of an implicational formula is determined from ones of its components at $v\succeq w$.
Let us apply the search procedure for ${\sf LKp}$ to a given sequent $S$
in which we analyze antecedent implicational formulas, but leave the succedent implicational formulas.
We obtain a deduction $Tr_{S}$, i.e., a tree of sequents which is locally correct with respect to inference rules.
If every leaf sequent in the tree $Tr_{S}$ is derivable, then so is the given sequent.
Here each leaf sequent is of the form $\Gam\Rarw\Del,\Del_{\supset}$
with the set $\Del_{\supset}$ of succedent implicational formulas and 
it is saturated in the sense of Definition \ref{df:saturated} below.
By saturation we see that 
there exists a formula $(\alp\supset\bet)\in\Del_{\supset}$ such that $\alp,\Gam\Rarw\bet$ is intuitionistically derivable
provided that $\Gam\Rarw\Del,\Del_{\supset}$ is intuitionistically derivable, 
$\Gam\Rarw\Del$ is not an axiom and $\Del_{\supset}\neq\emptyset$.
We keep on examining sequents $S_{\alp\supset\bet}=(\alp,\Gam\Rarw\bet)$ for each $(\alp\supset\bet)\in\Del_{\supset}$.
This yields deductions $Tr_{S_{\alp\supset\bet}}$, and so forth.
This means that we are constructing a tree $TR(S_{0})$ of deductions for a given sequent $S_{0}$.
Each node $\sig$ in the tree $TR(S_{0})$ corresponds to a deduction $Tr_{S(\sig)}$.
From the whole tree $TR(S_{0})$ one can extract either a derivation of $S_{0}$ or a Kripke model
in which $S_{0}$ is false, cf.\,Theorem \ref{th:intcomplete}.
Now details follow.
\\

Given a sequent $S=(\Gam\Rarw\Del)$, let $Tr_{S}$ denote the deduction of $S$ constructed 
in a bottom-up manner.

We analyze each formula in a sequent leaving implicational formulas in succedents. 
A formula is marked with a circle to indicate that the formula has not yet been analyzed.
No implicational formula in succedents is marked.
$\alp^{\circ}$ indicates that the formula $\alp$ has not been analyzed.
$\Gam^{\times}$ denotes the set of formulas obtained from formulas in $\Gam$ by erasing the circle.

\bdf\label{df:Sr}
{\rm Let $S=(\Gam\Rarw\Del)$ be a sequent.
$S_{\supset}$ denotes the set of implicational succedent formulas in $\Del$, while
 $S_{\supset}^{r}$ denotes the set of implicational formulas $\alp$ in $S_{\supset}$ such that
$\alp^{{\tt p}}\not\subset\Gam$ or $\alp^{{\tt c}}\not\in\Del$.
}
\edf
The condition 
$\alp^{{\tt p}}\subset\Gam \spand \alp^{{\tt c}}\in\Del$ means that the implicational succedent formula $\alp$ has 
been already analyzed
in the whole tree $TR(S_{0})$ of deductions defined in Definition \ref{df:TRp} below.

\bdf\label{df:saturated}
{\rm 
A sequent $\Gam\Rarw\Del$ is \textit{saturated} if the following five conditions are met:
\benu
\item
if $(\alp\lor\bet)\in\Gam$, then $\{\alp,\bet\}\cap\Gam^{\times}\neq\emptyset$,
\item
if $(\alp\lor\bet)\in\Del$, then $\{\alp,\bet\}\subset\Del^{\times}$,
\item
if $(\alp\land\bet)\in\Gam$, then $\{\alp,\bet\}\subset\Gam^{\times}$,
\item
if $(\alp\land\bet)\in\Del$, then $\{\alp,\bet\}\cap\Del^{\times}\neq\emptyset$,
and
\item
if an implicational formula $\alp$ is in $\Gam$, then $\alp^{{\tt p}}\subset\Del^{\times}$ or $\alp^{{\tt c}}\in\Gam^{\times}$.
\eenu

A saturated sequent $\Gam\Rarw\Del$ is \textit{fully analyzed} 
if it is not an axiom, i.e., $\Gam^{\times}\cap\Del^{\times}\cap Atm=\emptyset$ and $\bot\not\in\Gam^{\times}$,
and every marked formula $\alp^{\circ}$ in $\Gam\cup\Del$ is either an atom $p^{\circ}$
or $\bot^{\circ}$.
}
\edf

Note that these conditions for sequents to be saturated are only for unmarked formulas.


\bdf\label{df:inversionIp}
{\rm
Given a saturated sequent $S$, put the sequent $S$ at the root of the tree $Tr_{S}$, and
as far as one of the following inversion steps can be performed, continue it to construct the tree $Tr_{S}$.

Let us define inversion steps.
\benu
\item
If an antecedent $\Gam$ contains a conjunction $(\alp_{0}\land\alp_{1})^{\circ}$ with a circle,
then erase the circle and add marked conjuncts when these have not yet been analyzed: 
for $(\alp_{0}\land\alp_{1})^{\circ}\not\in\Gam_{1}$,
\[
\infer{(\alp_{0}\land\alp_{1})^{\circ},\Gam_{1}\Rarw\Del}
{\{\alp_{i}^{\circ}: \alp_{i}\not\in\Gam_{1}, i\in\{0,1\}\}\cup\Gam_{1}\cup\{\alp_{0}\land\alp_{1}\}\Rarw\Del}
\]
For simplicity let us denote it by
\[
\infer{(\alp_{0}\land\alp_{1})^{\circ},\Gam_{1}\Rarw\Del}
{\alp_{0}^{\st},\alp_{1}^{\st},\Gam_{1},\alp_{0}\land\alp_{1}\Rarw\Del}
\]
where $\alp_{i}^{\st}=\alp_{i}^{\circ}$ if $\alp_{i}\not\in\Gam_{1}$, and
$\alp_{i}^{\st}$ is absent else.
\item
If the succedent $\Del$ contains a conjunction $(\alp_{0}\land\alp_{1})^{\circ}$,
then erase the circle, add unmarked conjuncts and starred conclusions: 
for $(\alp_{0}\land\alp_{1})^{\circ}\not\in\Del_{1}$,
\[
\infer{\Gam\Rarw\Del_{1},(\alp_{0}\land\alp_{1})^{\circ}}
{
\Gam\Rarw \alp_{0}\land\alp_{1},\Del_{1},\alp_{0},\alp_{0}^{{\tt c}\st}
&
\Gam\Rarw \alp_{0}\land\alp_{1},\Del_{1},\alp_{1},\alp_{1}^{{\tt c}\st}
}
\]
where $\alp_{i}^{{\tt c}\st}\equiv(\alp_{i}^{{\tt c}})^{\circ}$ if $\alp_{i}^{{\tt c}}\not\in\Del_{1}$, and 
$\alp_{i}^{{\tt c}\st}$ is absent else.
\item
If an antecedent $\Gam$ contains a disjunction $(\alp_{0}\lor\alp_{1})^{\circ}$ with a circle,
then erase the circle and add starred disjuncts: for $(\alp_{0}\lor\alp_{1})^{\circ}\not\in\Gam_{1}$,
\[
\infer{(\alp_{0}\lor\alp_{1})^{\circ},\Gam_{1}\Rarw\Del}
{
\alp_{0}^{\st},\Gam_{1},\alp_{0}\lor\alp_{1}\Rarw\Del
&
\alp_{1}^{\st},\Gam_{1},\alp_{0}\lor\alp_{1}\Rarw\Del
}
\]
\item
If the succedent $\Del$ contains a disjunction $(\alp_{0}\lor\alp_{1})^{\circ}$ with a circle,
then erase the circle, add unmarked disjuncts and starred conclusions: 
for $(\alp_{0}\lor\alp_{1})^{\circ}\not\in\Del_{1}$,
\[
\infer{\Gam\Rarw\Del_{1},(\alp_{0}\lor\alp_{1})^{\circ}}
{
\Gam\Rarw \alp_{0}\lor\alp_{1},\Del_{1},\alp_{0},\alp_{1},\alp_{0}^{{\tt c}\st},\alp_{1}^{{\tt c}\st}
}
\]
\item
Suppose an antecedent $\Gam$ contains an implication $\alp^{\circ}$ with a circle.
Then erase the circle, add unmarked premisses and starred conclusions of premisses to the succedent,
and add starred conclusion to the antecedent:
for $\alp^{\circ}\not\in\Gam_{1}$,
\[
\infer{\alp^{\circ},\Gam_{1}\Rarw\Del}
{
\{
\alp,\Gam_{1}\Rarw\Del,\bet,\bet^{{\tt c}\st}
: \bet\in\alp^{{\tt p}}\}
&
\alp^{{\tt c}\st},\alp,\Gam_{1}\Rarw\Del
}
\]
\eenu
}
\edf
Observe that if the lower sequent enjoys the the following condition (\ref{eq:saturated}),
then so do the upper sequents in each inversion step.
\beqn\label{eq:saturated}
\mbox{for any implicational formula $\alp$ in the succedent $\Del$, $\alp^{{\tt c}}\in\Del^{\times}$.}
\eeqn

\bdf
{\rm
Let us define numbers $b(\alp)$ for marked or unmarked formulas $\alp$ recursively as follows.
\benu
\item
$b(\alp)=b(p^{\circ})=b(\bot^{\circ})=0$ for any unmarked formulas $\alp$.
\item
$b((\alp\lor\bet)^{\circ})=b((\alp\land\bet)^{\circ})=b(\alp^{\circ})+b(\bet^{\circ})+1$.
\item
$b(\alp^{\circ})=\max(\{b(\bet^{{\tt c}\circ}):\bet\in\alp^{{\tt p}}\}\cup\{b(\alp^{{\tt c}\circ})\})+1$
for implicational formulas $\alp$.
\eenu
For sequents $S=(\Gam\Rarw\Del)$,
let $b(S)=\sum\{b(\alp):\alp\in\Gam\cup\Del\}$.
}
\edf

\bprp\label{prp:Trterminate}
The number decreases when we go up in the tree $Tr_{S}$.
Namely in the above inversion steps
\[
\infer{S_{0}}
{
S_{00}
&
(S_{01})
}
\]
$b(S_{0i})<b(S_{0})$ holds for $i\in\{0,1\}$.
\eprp

Hence the process terminates, and $Tr_{S}$ is a finite tree.
Assume that $S$ is a saturated sequent.
Then each leaf in $Tr_{S}$ is a saturated sequent, which is
an axiom, or a fully analyzed sequent.

A fully analyzed sequent $S=(\Gam,\Gam_{\supset}\Rarw\Del,\Del_{\supset}^{r})$ 
with $\Del_{\supset}^{r}=S_{\supset}^{r}\neq\emptyset$
is extended by an `inference rule' \textit{branching} $(\mbox{br})$,
where $\Gam_{\supset}$ is the set of antecedent unmarked implicational formulas.

$\Gam_{\supset}^{\circ}$ is obtained from $\Gam_{\supset}$ by marking each formula
to analyze these again,
$\Gam_{\supset}^{\circ}=\{\alp^{\circ}:\alp\in\Gam_{\supset}\}$.
Let $\alp^{{\tt p}\st}=\{\bet^{\circ}:\bet\in\alp^{{\tt p}}, \bet\not\in\Gam\cup\Gam_{\supset}\}$.
\[
\infer[(\mbox{br})]{\Gam,\Gam_{\supset}\Rarw\Del,\Del_{\supset}^{r}}
{
\{\alp^{{\tt p}\st}\cup\Gam\cup\Gam_{\supset}^{\circ}\Rarw\alp^{{\tt c}\circ}: \alp\in\Del_{\supset}^{r}\}
}
\]
where each formula $\alp\in\Del_{\supset}^{r}$ is said to be a \textit{major formula} of the rule $(\mbox{br})$.

This inference rule is a disjunctive one since if one of upper sequents is derivable, 
then so is the lower sequent using the inference rule $(R\supset)$.
Note that each upper sequent $S_{\alp}=(\alp^{{\tt p}\st}\cup\Gam\cup\Gam_{\supset}^{\circ}\Rarw\alp^{{\tt c}\circ})$ enjoys the condition (\ref{eq:saturated}) vacuously
since $(S_{\alp})_{\supset}=(S_{\alp})_{\supset}^{r}=\emptyset$.

\subsection{Construction of the tree of deductions}\label{subsec:IpTR}
Let us define a tree $TR(S_{0})$ of deductions for saturated sequents $S_{0}$.
Each node in the tree $TR(S_{0})$ corresponds to a deduction $Tr_{S}$.

${}^{<\ome}\ome$ denotes the set of finite sequences of natural numbers $\sig,\tau,\ldots$.
The empty sequence is denoted by $\emptyset$, and $(k_{0},\ldots,k_{n-1})*(j)=(k_{0},\ldots,k_{n-1},j)$.

\bdf\label{df:TRp}
{\rm Given a saturated sequent $S_{0}=(\Gam_{0}\Rarw\Del_{0})$, let us define a tree $TR(S_{0})\subset{}^{<\ome}\ome$, 
 and a labeling function $(S(\sig), g(\sig),d(\sig))$ for $\sig\in TR(S_{0})$,
where $S(\sig)$ is a saturated sequent, $g(\sig)\in\{\lor,\land,0,1\}$ is a (logic) gate and $d(\sig)$ is a deduction possibly with the branching rule
such that
\benu
\item
for each leaf $\sig$ in $TR(S_{0})$, $d(\sig)$ consists solely of the saturated sequent $S(\sig)$, and either 
 \benu
 \item
$S(\sig)$ is an axiom and $g(\sig)=1$ indicating that $S(\sig)$ is derivable, or 
 \item
 $S(\sig)$ is fully analyzed and $g(\sig)=0$ indicating that $S(\sig)$ is underivable, and
 \eenu
\item
for each internal node $\sig$ in $TR(S_{0})$, $g(\sig)\in\{\lor,\land\}$ and
 \benu
 \item
 if $g(\sig)=\land$, then 
 $S(\sig)$ is neither an axiom nor fully analyzed one, and 
 $d(\sig)=Tr_{S(\sig)}$, where $Tr_{S(\sig)}$ is the deduction for the sequent $S(\sig)$ defined in Definition \ref{df:inversionIp}.
 \item
  if $g(\sig)=\lor$, then $S(\sig)$ is a fully analyzed sequent such that $(S(\sig))_{\supset}^{r}\neq\emptyset$,
  and $d(\sig)$ is a deduction with a single inference rule $(\mbox{br})$
  with its lower sequent $S(\sig)$.
  \eenu
\eenu
Thus $TR(S_{0})$ is a $(\land,\lor)$-tree, and it is constructed inductively according to 
$\land${\bf -stage} or to $\lor${\bf -stage} below.
\\

\noindent
{\bf initial}.
First the empty sequence $\emptyset\in TR(S_{0})$ and $S(\emptyset)=S_{0}$.
If $S_{0}$ is an axiom, then $g(\emptyset)=1$.
If $S_{0}$ is fully analyzed with $(S_{0})^{r}_{\supset}=\emptyset$, then $g(\emptyset)=0$.
If $g(\emptyset)\in\{0,1\}$, then $d(\emptyset)$ is the deduction consisting solely of $S_{0}$.
If $S_{0}$ is fully analyzed with $(S_{0})_{\supset}^{r}\neq\emptyset$, then $g(\emptyset)=\lor$
and the tree is extended according to the $\lor${\bf -stage}.
Otherwise $g(\emptyset)=\land$ and the tree is extended according to the $\land${\bf -stage}.
\\

\noindent
$\land${\bf -stage}.
Suppose $\sig\in TR(S_{0})$ and $g(\sig)=\land$. Let $S(\sig)=S$.

Let $\{S_{i}\}_{i<I}\, (I>0)$ be an enumeration of all leaves in $d(\sig)=Tr_{S}$.
For each $i<I$, let $\sig*(i)\in TR(S_{0})$ with $S(\sig*(i))=S_{i}$.
If $S_{i}$ is an axiom, then $g(\sig*(i))=1$.
Otherwise $S_{i}$ is fully analyzed.
If $(S_{i})_{\supset}^{r}=\emptyset$, then $g(\sig*(i))=0$.
Otherwise let $g(\sig*(i))=\lor$ and the tree is extended according to the $\lor${\bf -stage}.
\\

\noindent
$\lor${\bf -stage}.
Suppose $\sig\in TR(S_{0})$ and $g(\sig)=\lor$. 

Let $S(\sig)=S$ be a fully analyzed sequent
$\Gam,\Gam_{\supset}\Rarw\Del,\{\alp_{j}\}_{j<J}\,(J>0)$
where $\Gam_{\supset}$ is the set of antecedent unmarked implicational formulas, 
and $\{\alp_{j}\}_{j<J}=S_{\supset}^{r}$.
Then $\sig*(j)\in TR(S_{0})$ for each $j<J$.
Also let $S(\sig*(j))=(\alp_{j}^{{\tt p}\st}\cup\Gam\cup\Gam_{\supset}^{\circ}\Rarw\alp_{j}^{{\tt c}\circ})$,
where
each unmarked and implicational formula $\alp\in\Gam_{\supset}$ is marked in $\Gam^{\circ}_{\supset}$
to be analyzed again.
$d(\sig)$ denotes the following deduction:
\[
\infer[(\mbox{br})]{\Gam,\Gam_{\supset}\Rarw\Del,\{\alp_{j}\}_{j<J}}
{
\{\alp_{j}^{{\tt p}\st}\cup\Gam\cup\Gam_{\supset}^{\circ}\Rarw\alp_{j}^{{\tt c}\circ}\}_{j<J}
}
\]

Let $g(\sig*(j))=1$ if $S(\sig*(j))$ is an axiom, and $g(\sig*(j))=0$ if it is fully analyzed with 
$(S(\sig*(j)))_{\supset}^{r}=\emptyset$.
Otherwise let $g(\sig*(j))=\land$ and the tree is extended according to the $\land${\bf -stage}.
}
\edf

For a formula $\alp$, let $Sbfml_{\supset}^{+}(\alp)$ [$Sbfml_{\supset}^{-}(\alp)$] denote a set of positive 
[negative] implicational subformulas in $\alp$ defined by simultaneous recursion as follows.
\benu
\item
$Sbfml_{\supset}^{\pm}(p)=Sbfml_{\supset}^{\pm}(\bot)=\emptyset$.
\item
$Sbfml_{\supset}^{\pm}(\alp_{0}\lor\alp_{1})=Sbfml_{\supset}^{\pm}(\alp_{0}\land\alp_{1})=Sbfml_{\supset}^{\pm}(\alp_{0})\cup Sbfml_{\supset}^{\pm}(\alp_{1})$.
\item
Let $\alp$ be an implicational formula.
\\
$Sbfml_{\supset}^{+}(\alp)=\bigcup\{Sbfml_{\supset}^{-}(\bet): \bet\in\alp^{{\tt p}}\}\cup Sbfml_{\supset}^{+}(\alp^{{\tt c}})\cup\{\alp\}$.
\\
$Sbfml_{\supset}^{-}(\alp)=\bigcup\{Sbfml_{\supset}^{+}(\bet): \bet\in\alp^{{\tt p}}\}\cup Sbfml_{\supset}^{-}(\alp^{{\tt c}})$.
\eenu
For sequents $S=(\Gam\Rarw\Del)$,
let $Sbfml_{\supset}^{+}(S)=\bigcup\{Sbfml_{\supset}^{-}(\alp):\alp\in\Gam\}\cup\bigcup\{Sbfml_{\supset}^{+}(\bet):\bet\in\Del\}$.

\blem\label{prp:termination}
The whole process generating the tree $TR(S_{0})$ terminates, and
the number of $\lor$-gates along any branch in the tree $TR(S_{0})$ is at most 
the cardinality $\#Sbfml_{\supset}^{+}(S_{0})$ of the set $Sbfml_{\supset}^{+}(S_{0})$.
\elem
\bprf
Consider a branching rule in the tree $TR(S_{0})$.
\[
\infer[(\mbox{br})]{\Gam,\Gam_{\supset}\Rarw\Del,\{\alp_{j}\}_{j<J}}
{
\cdots
&
\infer*[\pi]{\alp_{j}^{{\tt p}\st}\cup\Gam\cup\Gam_{\supset}^{\circ}\Rarw\alp_{j}^{{\tt c}\circ}}
 {
  \infer[(\mbox{br})\, I]{S}{}
  }
&
\cdots
}
\]
Each major formula $\alp_{j}$ of the branching rule is in the set $Sbfml^{+}_{\supset}(S_{0})$.
In the branch $\pi$ up to an upper sequent 
$S_{j}=(\alp_{j}^{{\tt p}\st}\cup\Gam\cup\Gam_{\supset}^{\circ}\Rarw\alp_{j}^{{\tt c}\circ})$,
each formula in the set $\alp_{j}^{{\tt p}}$ is in the antecedents of sequents $S$ in $\pi$.
Moreover $S$ enjoys the condition (\ref{eq:saturated}) since $S_{j}$ enjoys it and
it is preserved upward.

Suppose that there occurs a branching rule $I$ on $\pi$ above the upper sequent $S_{j}$ one of whose major formulas
is the formula $\alp_{j}$.
Let $S$ be the lower sequent of the branching rule $I$.
Then the formula  $\alp_{j}$ is in the succedent of $S$, 
and hence $\alp_{j}^{{\tt c}}$ is in the succedent by  (\ref{eq:saturated}).
Furthermore
each formula in the set $\alp_{j}^{{\tt p}}$ is in the antecedent of sequent $S$.
Hence $\alp_{j}\not\in S_{\supset}^{r}$,
and $\alp_{j}$ is not a major formula of the branching rule with the lower sequent $S$.

Therefore there are at most $\#Sbfml_{\supset}^{+}(S_{0})$ applications of the branching rules along a branch in $TR(S_{0})$.
This means that the number of $\lor$-gates along any branch in the tree $TR(S_{0})$ is at most $\#Sbfml_{\supset}^{+}(S_{0})$.
Since each deduction $Tr_{S}$ is finite, the whole process generating the tree $TR(S_{0})$ terminates.
\eprf
\\

\noindent
Let us compute the value of the $(\land,\lor)$-tree $TR(S_{0})$ with gates $g(\sig)$.
Let $v(\sig)$ denote the value of $\sig\in TR(S_{0})$.
If the value $v(\emptyset)$ is $1$, then $S_{0}$ is derivable where a derivation of $S_{0}$ is obtained by
putting the deduction $Tr_{S(\sig)}$ of $S(\sig)$ from $\{S_{i}\}_{i<I}$ 
\[
\begin{array}{ccc}
\cdots & S_{i} & \cdots
\\
\searrow & \downarrow & \swarrow
\\
& (S(\sig),\land) &
\end{array}
\]
to the $\land$-node $\sig$, and choosing one of upper sequents $S(\sig*(i,j_{i}))$ such that
$v(\sig*(i,j_{i}))=1$ for each
lower sequent $S(\sig*(i))$ of $(\mbox{br})$, i.e.,
$g(\sig*(i))=\lor$.

In what follows consider the case when the value $v(\emptyset)$ of $TR(S_{0})$ is $0$.
In a bottom-up manner let us shrink the tree $TR(S_{0})$ to a tree $T\subset TR(S_{0})$ as follows.
Simultaneously a set $V_{T}(\sig)$ of atoms is assigned.
For each node $\sig\in T$ $g(\sig)\neq\lor$.
First $\emptyset\in T$.

Suppose $v(\sig)=0$ for a node $\sig\in T$ with $g(\sig)=\land$.
Pick a son $\sig*(i)$ such that $v(\sig*(i))=0$, and identify the node $\sig*(i)$ with $\sig$.
This means that we have chosen an underivable sequent $S(\sig*(i))$, which is a non-axiom leaf
in the deduction $Tr_{S(\sig)}$.
Let $V_{T}(\sig)=\Gam(\sig)^{\times}\cap Atm$ where $\Gam(\sig)\Rarw\Del(\sig)$ denotes the leaf sequent
$S(\sig*(i))$ chosen from $Tr_{S(\sig)}$.
If $g(\sig*(i))=0$, then $\sig$ will be a leaf in $T$.
Otherwise $g(\sig*(i))=\lor$, and keep its sons in a shrunken tree, i.e., $\sig*(i,j)\in T$.

Thus we have defined a Kripke model $\langle T,\subset_{e},V_{T}\rangle$ where 
$\sig\subset_{e} \tau$ iff $\sig$ is an initial segment of $\tau$.
In this case we say that $\tau$ is an extension of $\sig$.
By definition, $\sig\subset_{e}\sig$.

For each $\sig\in T$ with $g(\sig)=\land$, $\Gam(\sig)\Rarw\Del(\sig)$ denotes the leaf sequent
$S(\sig*(i))$ chosen from $Tr_{S(\sig)}$.

From the construction we see readily the followings for any $\sig,\tau\in T$.
\benu
\item
$\sig\subset_{e}\tau\Rarw \Gam(\sig)^{\times}\subset\Gam(\tau)^{\times}$.
\item
$\Gam(\sig)\Rarw\Del(\sig)$ is saturated.
\item
 if an implicational formula $\alp\in \Del(\sig)$, then there exists an extension $\tau\in T$ of $\sig$ such that
 $\alp^{{\tt p}}\subset \Gam(\tau)^{\times}$ and $\alp^{{\tt c}}\in\Del(\tau)^{\times}$.

\item
$\Gam(\sig)^{\times}\cap\Del(\sig)^{\times}$ has no common atom, and $\bot\not\in\Gam(\sig)^{\times}$.
\eenu

\bprp
If $\alp\in\Gam(\sig)^{\times}$ $[\alp\in\Del(\sig)^{\times}]$, then $\sig\models\alp$ $[\sig\not\models\alp]$, resp. in the
 Kripke model $\langle T,\subset_{e},V_{T}\rangle$.
Hence $\sig\models\bigwedge\Gam(\sig)$ and $\sig\not\models\bigvee\Del(\sig)$.
\eprp
\bprf
By simultaneous induction on $\alp$ using the above facts.
\eprf

\bth\label{th:intcomplete}{\rm (Sch\"utte's dichotomy)}\\
For every saturated sequent $S_{0}$,
$v(\emptyset)=1$ iff ${\sf LJpm}^{{\tt c}}\vdash S_{0}$.

Specifically if $v(\emptyset)=0$, 
then each Kripke model $\langle T,\subset_{e},V_{T}\rangle$ falsifies the given saturated sequent $S_{0}$,
no matter which non-axiom leaves are chosen from $Tr_{S_{i}}$.

On the contrary, if $v(\emptyset)=1$, then
we can extract a (cut-free) derivation of $S_{0}$ by
choosing a derivable sequent from each $(\mbox{br})$.

Hence ${\sf LJpm}^{{\tt c}}$ is intuitionistically complete in the sense that
any intuitionistically valid sequent is derivable in ${\sf LJpm}^{{\tt c}}$.

Moreover ${\sf LJpm}^{{\tt c}}$ admits the Hauptsatz, i.e.,
the cut rule is admiissible:
if both of the sequents $\Gam_{0}\Rarw\Del_{0},\alp$ and $\alp,\Gam_{1}\Rarw\Del_{1}$ 
are derivable in ${\sf LJpm}^{{\tt c}}$,
then so is the sequent $\Gam_{0},\Gam_{1}\Rarw\Del_{0},\Del_{1}$.
\end{theorem}
\bprf
Let $S=(\Gam\Rarw\Del,\Del_{\supset})$ be a sequent without circles,
where $\Del_{\supset}=S_{\supset}$ is the set of implicational formulas in the succedent.
Then $S_{0}=(\Gam^{\circ}\Rarw\Del^{\circ},\Del_{\supset})$ is saturated.
\eprf

\bdf\label{df:De}
{\rm For a sequent $S_{0}=(\Gam_{0}\Rarw \Del_{0})$, let 
$S_{0}^{\circ}=(\Gam_{0}^{\circ}\Rarw \Del_{1}^{\circ},\Del_{\supset})$ be the saturated sequent,
where $\Del_{0}=\Del_{1}\cup\Del_{\supset}$ with $\Del_{\supset}=(S_{0})_{\supset}$.
Then $De(S_{0})$ denotes the deduction with the branching rule $(\mbox{br})$, obtained from
 $TR(S_{0}^{\circ})$ by carrying out the intermediate deductions $Tr_{S(\sig)}$ for 
$\sig\in TR(S_{0}^{\circ})$ with $g(\sig)=\land$.
}
\edf

\bcor\label{prp:depthIp}
Let $S_{0}$ be a sequent and $\pi$ be a branch in the tree $De(S_{0})$.
\benu
\item
For each implicational formula $\alp$ in $Sbfml^{+}_{\supset}(S_{0}^{\circ})$,
there is at most one application of a right rule $(R\supset)$ with the major formula $\alp$ along the branch $\pi$.
Hence the number of applications of the rule $(R\supset)$ along the branch $\pi$ is at most 
$\# Sbfml^{+}_{\supset}(S_{0}^{\circ})$.

\item
The number of applications of the left rule $(L\supset)$ along the branch $\pi$ is at most 
$1+\# Sbfml^{+}_{\supset}(S_{0}^{\circ})$.
\item
The depth of the tree $De(S_{0})$ is bounded by
$n(S_{0})\cdot (1+\# Sbfml^{+}_{\supset}(S_{0}^{\circ}))$, where
$n(S_{0})$ denotes the total number of occurrences of connectives $\lor,\land,\supset$
in $S_{0}$.
\eenu
Each intuitionistically valid sequent $S_{0}$ has a derivation in ${\sf LJpm}^{{\tt c}}$ enjoying the
above conditions with respect to the number of applications of rules $(R\supset),(L\supset)$ along any branch and the depth of the derivation.
\ecor
\bprf
These are seen from the proofs of Proposition \ref{prp:Trterminate} and Lemma \ref{prp:termination} using the fact
$b(S)\leq n(S)$.
\eprf

\bcor\label{cor:pedagogic}
Let $S_{0}=(\Gam_{0}\Rarw\Del_{0})$ be a sequent.
Let $\langle W,\preceq,V\rangle$ be a Kripke model falsifying the sequent $S_{0}$.
Let $w_{0}\in W$ be such that
$w_{0}\models\bigwedge\Gam$ and $w_{0}\not\models\bigvee\Del$.

Then we can choose a Kripke model $\langle T,\subset_{e},V_{T}\rangle$ for which
there exists an $h:T\to W$ such that $\sig\subset_{e}\tau\Rarw w_{0}\preceq h(\sig)\preceq h(\tau)$
and $V_{T}(\sig)\subset V(h(\sig))$ for any $\sig,\tau\in T$.
\ecor
\bprf
Choose a tree $T$, i.e., pick leaves of deductions bottom-up manner in $De(S_{0})$ 
according to the given Kripke model $\langle W,\preceq,V\rangle$.
For example if $w\not\models(\alp_{0}\land\alp_{1})$, then pick an $i$ such that $w\not\models\alp_{i}$, and
go to the $i$-th branch.
Suppose $w\models(\alp\supset\bet)$.
If $w\not\models\alp$, then go to the left. 
Otherwise we have $w\models\bet$, and go to the right.

When we reach a leaf $\sig*(i)$ of $Tr_{S(\sig)}$ with $g(\sig*(i))=\lor$ and $\alp\equiv(\gam\supset\del)$ 
is one of the succedent formula in $(S(\sig*(i)))_{\supset}^{r}$,
supposing $w\not\models\alp$, and let $v\succeq w$ be such that $v\models\gam$ and $v\not\models\del$.
Put $h(\sig*(i,j))=v$ where in $\sig*(i,j)$, $\alp$ is analyzed.

When we reach a leaf $S(\sig*(i))=(\Gam(\sig)\Rarw\Del(\sig))$ of $Tr_{S(\sig)}$ with $g(\sig*(i))=0$ and $w=h(\sig)$, then 
$w\models p$ with $p^{\circ}\in\Gam(\sig)$.
This means that any $p\in V_{T}(\sig)$ is in $V(w)=V(h(\sig))$, and hence $V_{T}(\sig)\subset V(h(\sig))$.
\eprf

\subsection{The intuitionistic propositional logic is in PSPACE}\label{subsec:IpPSPACE}

In this subsection we describe a PSPACE-algorithm deciding the deducibility of the given sequent $S_{0}$ in 
${\sf LJpm}^{{\tt c}}$.
This is a result due to R. Ladner\cite{Ladner}, who shows that a modal logic {\sf S4} is in PSPACE, and it is well known that
{\sf S4} interprets the intuitionistic propositional logic  {\sf Ip} linearly.
J. Hudelmaier\cite{Hudelmaier} sharpens the result by giving an $O(n\log n)$-space decision procedure for {\sf Ip}.

\bcor\label{cor:PSPACEIp}{\rm (Ladner\cite{Ladner})}\\
The intuitionistic propositional logic  {\sf Ip} is in PSPACE.
\ecor

Recall that $De(S_{0})$ denotes the deduction obtained from $TR(S_{0}^{\circ})$ by carrying out the
 intermediate deductions $Tr_{S(\sig)}$ for 
$\sig\in TR(S_{0}^{\circ})$ with $g(\sig)=\land$.
Let $\# S$ be the \textit{size} of the sequents $S$, which is the total number of occurrences of symbols in $S$.

\bprp\label{prp:maxsize}
Let $\vec{S}=S_{0},S_{1},\ldots, S_{n-1}$ be a branch in $De(S_{0})$, where $S_{i+1}$ is an upper sequent of an inference rule
with its lower sequent $S_{i}$.
Then $\#\vec{S}:=\sum_{i< n}\# S_{i}$ is bounded by a (quartic) polynomial in the size $\# S_{0}$ of the given sequent $S_{0}$.
\eprp
\bprf
From Corollary \ref{prp:depthIp} and $n(S_{0})\leq \# S_{0}$
we see that the length $n$ of branches $\vec{S}$ is at most $(1+\# Sbfml_{\supset}^{+}(S_{0}))\cdot \# S_{0}$,
which is
bounded by a quadratic polynomial in $\# S_{0}$.
On the other side the maximal size of sizes $\# S_{i}$ of sequents $S_{i}$ is bounded by a quadratic polynomial, too,
since each $S_{i}$ is essentially a sequence of subformulas of formulas in $S_{0}$.
\eprf
\\

\noindent
Let us traverse sequents in the tree $De(S_{0})$ starting from the root $S_{0}$ as follows.
Let $S$ be the current sequent.
\bdes
\item[Case 1]
If $S$ is a lower sequent of an inference rule in $De(S_{0})$, then visit the leftmost upper sequent next.

\item[Case 2]
Otherwise $S$ is a leaf in $De(S_{0})$, and in $TR(S_{0})$.
Let $\sig\in TR(S_{0})$ be the node such that $S(\sig)=S$.
$S(\sig)$ is a leaf in a deduction $Tr_{S(\tau)}$ for a $\tau\in TR(S_{0})$
with $g(\tau)=\land$.

 \item[Case 2.1]
First consider the case when $g(\sig)=1$.
\item[Case 2.1.1]
Let $T$ be the uppermost sequent below $S$ in $Tr_{S(\tau)}$
such that there is an upper sequent $S^{\prime}$ of the inference rule with its lower sequent $T$, which we have not  yet visited.
Next let us visit the leftmost such sequent $S^{\prime}$ if such a sequent exists.
\[
Tr_{S(\tau)}=
\infer*{}
{
 \infer{T}
 {
 \cdots
 &
 \infer*{}{S}
 &
 S^{\prime}
 &
 \cdots
 }
}
\]
\item[Case 2.1.2]
Suppose that there is no such sequent.
This means that $S(\tau)$ is derivable.
If $\tau=\emptyset$, then we are done.
Otherwise $S(\tau)$ is an upper sequent of a $(\mbox{br})$, and we see that the lower sequent $S(\rho)$ of the $(\mbox{br})$
is derivable.
Let us change the gate $g(\rho)=\lor$ to $g(\rho)=1$, and continue the search for the next visiting sequent
in the deduction $Tr_{S(\kap)}$, where $S(\rho)$ is a leaf in the deduction $Tr_{S(\kap)}$.

\item[Case 2.2]
Second consider the case when $g(\sig)=0$.
This means that $S(\tau)$ is underivable.
If $\tau=\emptyset$, then we are done.
Otherwise $S(\tau)$ is an upper sequent of a $(\mbox{br})$.
\item[Case 2.2.1]
If $S(\tau)$ is not the rightmost upper sequent, then visit the next right one.
\item[Case 2.2.2]
Otherwise the lower sequent $S(\rho)$ of the $(\mbox{br})$ is underivable.
Let us change the gate $g(\rho)=\lor$ to $g(\rho)=0$, and continue the search for the next visiting sequent
in the deduction $Tr_{S(\kap)}$, where $S(\rho)$ is a leaf in the deduction $Tr_{S(\kap)}$.

\edes

In the PSPACE-algorithm,
we record sequences $\vec{S}=S_{0},S_{1},\ldots,S_{n-1}$ of sequents on a tape, where
$\vec{S}$ is an initial segment of a branch in $De(S_{0})$.
The next sequence $\vec{S}^{\prime}$ is recursively computed as follows.
If the tail $S_{n-1}$ is a lower sequent of an inference rule in $De(S_{0})$, then 
$\vec{S}^{\prime}=\vec{S}*(S_{n})$, i.e., extend the sequence $\vec{S}$
by adding the leftmost upper sequent $S_{n}$ as a tail, cf.\,{\bf Case 1} in the traversal.

Suppose $S_{n-1}=S(\sig)$ is a leaf in a deduction $Tr_{S(\tau)}$.

First consider the case when $g(\sig)=1$.
If there is an $S_{i}\, (i<n-1)$ such that $S_{i+1}$ is not the rightmost upper sequent,
then break the sequence $\vec{S}$ at $S_{i}$ and put the next right upper sequent $S^{\prime}$,
$\vec{S}^{\prime}=S_{0},\ldots,S_{i},S^{\prime}$ for the maximal such $i$, cf.\,{\bf Case 2.1.1}.
\[
\infer{S_{i}}
{
\cdots
&
S_{i+1}
&
S^{\prime}
&
\cdots
}
\]
Suppose there is no such $S_{i}$.
If $\tau=\emptyset$, then halt and print `DERIVABLE'.
Otherwise $S(\tau)$ is an upper sequent of a $(\mbox{br})$ with the lower sequent $S(\rho)$.
Continue the computation of the next sequence $\vec{S}^{\prime}$ 
in the deduction $Tr_{S(\kap)}$, where $S(\rho)$ is a leaf in the deduction $Tr_{S(\kap)}$, cf.\,{\bf Case 2.1.2}.

Second consider the case when $g(\sig)=0$.
If $\tau=\emptyset$, then halt and print `UNDERIVABLE'.
Otherwise $S(\tau)$ is an upper sequent of a $(\mbox{br})$.
If $S(\tau)$ is not the rightmost upper sequent, then 
break the sequence $\vec{S}$ at $S_{i}$ and put the next right upper sequent $S^{\prime}$,
$\vec{S}^{\prime}=S_{0},\ldots,S_{i},S^{\prime}$, where
$S_{i+1}=S(\tau)$ and $S_{i}$ is the lower sequent of the $(\mbox{br})$,
cf.\,{\bf Case 2.2.1}.
\[
\infer[(\mbox{br})]{S_{i}}
{
\cdots
&
S_{i+1}(=S(\tau))
&
S^{\prime}
&
\cdots
}
\]
Otherwise continue the computation of the next sequence $\vec{S}^{\prime}$ 
in the deduction $Tr_{S(\kap)}$, 
where the lower sequent $S_{i}$ of the $(\mbox{br})$ is a leaf in the deduction $Tr_{S(\kap)}$.

In each moment, we see from Proposition \ref{prp:maxsize} that
the size $\#\vec{S}$ of the recorded sequence $\vec{S}$ is bounded by 
a (quartic) polynomial in the size $\# S_{0}$ of the given sequent $S_{0}$.
Therefore the algorithm runs in a polynomial space in the size $\# S_{0}$ of the input $S_{0}$.

\section{Decidable fragments of {\sf Iq}}\label{sec:Iq}

A \textit{relational} language $\calL$ of the predicate logic (without equality) consists of
\textit{propositional connectives} $\bot,\lor,\land,\supset$,
\textit{quantifiers} $\exi,\fal$, (finite) list of \textit{predicate symbols} $R,\ldots$, \textit{individual contants} $c,\ldots$
\textit{Free variables} $a_{i}\, (i\in\ome)$ are denoted by $a,\ldots$, and \textit{bound variables} $x_{i}, (i\in\ome)$
are denoted by $x,y,z$.
$FV=\{a_{i}:i\in\ome\}$ denotes the set of free variables.
In this section $Atm$ denotes the set of atomic formulas $R(t_{1},\ldots,t_{n})$, where $t_{i}$ is a term, i.e., 
either an individual constant or a free variable.
A formula is said to be \textit{relational} if it is a formula in a relational language.
 {\sf Iq} denotes the intuitionistic predicate logic over a relational language.

\subsection{Sequent calculus ${\sf LJm}^{{\tt c}}$ for {\sf Iq}}\label{subsec:LJm}
${\sf LJm}^{{\tt c}}$ is obtained from ${\sf LJpm}^{{\tt c}}$  by adding the following inference rules for quantifiers
with their \textit{eigenvariables} $a$:
\[
\infer[(L\exi)]{\Gam\Rarw\Del}
{\alp(a),\Gam\Rarw\Del}
\]
where $\exi x\,\alp(x)\in\Gam$ and the free variable $a$ does not occur in $\Gam\cup\Del$.

\[
\infer[(R\exi)]{\Gam\Rarw\Del}
{
\Gam\Rarw\Del\cup\{\alp(t_{i}),\alp(t_{i})^{{\tt c}}\}_{i}
}
\mbox{with } \exi x\,\alp(x)\in\Del \mbox{ and a non-empty list of terms } \{t_{i}\}_{i}
\]
\[
\infer[(L\fal)]{\Gam\Rarw\Del}
{\{\alp(t_{i})\}_{i}\cup\Gam\Rarw\Del}
\mbox{with } \fal x\,\alp(x)\in\Gam  \mbox{ and a non-empty list of terms } \{t_{i}\}_{i}
\]
\[
\infer[(R\fal)]{\Gam\Rarw\Del}
{
\Gam\Rarw\alp(a),\alp(a)^{{\tt c}}
}
\]
where $\fal x\,\alp(x)\in\Del$ and the free variable $a$ does not occur in $\Gam\cup\{\fal x\,\alp(x)\}$.

A \textit{Kripke model} for a relational language $\calL$ is a quadruple $\langle W,\preceq,D,I\rangle$, where 
$\langle W, \preceq\rangle$ is a Kripke frame.
This has to enjoy the following for $w\preceq v$.
$D:W\to \calP(X)$ for a set $X$ such that $\emptyset\neq D(w)\subset D(v)$, and 
for each $w\in W$,
$I(w)$ is an $\calL$-structure with the universe $D(w)$ and relations $R^{w}\subset D(w)^{n}$ for $n$-ary predicate symbols $R\in\calL$
and elements $c^{w}\in D(w)$ for individual contants $c\in\calL$ 
such that
$R^{w}\subset R^{v}$,
and $c^{w}=c^{v}$.

For closed formulas $\alp\in \calL(X)$ and $w\in W$, $w\models\alp$ is defined recursively.

\benu
\item
$w\models R(c_{1},\ldots,c_{n})$ iff $(c_{1},\ldots,c_{n})\in R^{w}$.
\item
$w\models\alp\lor\bet$ iff $w\models\alp$ or $w\models\bet$.
\item
$w\models\alp\land\bet$ iff $w\models\alp$ and $w\models\bet$.
\item
$w\not\models\bot$.
\item
$w\models(\alp\supset\bet)$ iff for any $v\succeq w$, if $v\models\alp$, then $v\models\bet$.
\item
$w\models \exi x\,\alp(x)$ iff there exists a $c\in D(w)$ such that $w\models\alp(c)$.
\item
$w\models\fal x\,\alp(x)$ iff for any $v\succeq w$ and any $c\in D(v)$, $v\models\alp(c)$.
\eenu

\subsection{Proof searches for fragments}\label{subsec:Iq+PSPACE}

A formula is said to be \textit{positive} (with respect to quantifiers) iff
any universal quantifier [existential quantifier] occurs only positively [negatively] in it, resp.
Here positive/negative occurrence of quantifiers is meant in the usual classical sense.

G. Mints\cite{MintsSol} showed that it is decidable whether or not a given positive formula is intuitionisitically derivable,
cf.\,\cite{Dowek} for an alternative proof.

Let us introduce classes $\fal^{(++)}\exi^{(-)}$ of formulas and its subclass $\fal^{(++,-)}\exi^{(+)}$.

\bdf
{\rm
\benu
\item
A relational formula is defined to be in the class $\fal^{(++,-)}\exi^{(+)}$
iff each positively occurring universal quantifier occurs strictly positive, and existential quantifiers occur only positively in it.

Universal quantifiers may occur negatively.
\item
A relational formula is defined to be in the class $\fal^{(++)}\exi^{(-)}$
iff universal quantifiers occur only strictly positive, and existential quantifiers occur only negatively in it.
\eenu
}
\edf
For example
a positive formula $(\fal x(P(x)\to q)\to q) \to q$ is not in these classes,  but equivalent to the formula
$((\exi x\, P(x)\to q)\to q) \to q$ in the $\fal^{(++)}\exi^{(-)}$.

In this subsection we show that 
there is an algorithm running in polynomial space, which decides the intuitionistic derivability of
formulas in the class $\fal^{(++)}\exi^{(-)}$,
and there is an algorithm running in exponential space, which decides the intuitionistic derivability of
formulas in the class 
$\fal^{(++,-)}\exi^{(+)}$.

A sequent is said to be in the class $\fal^{(++,-)}\exi^{(+)}$ iff any succedent formula in it is 
in $\fal^{(++,-)}\exi^{(+)}$, 
universal quantifiers occur only positively and
existential quantifiers occur only negatively in its antecedent formulas.
A sequent is said to be in the class $\fal^{(++)}\exi^{(-)}$ iff any succedent formula in it is 
in $\fal^{(++)}\exi^{(-)}$, and in its antecedent formulas
no universal quantifier occurs, and existential quantifiers occur only positively.

A formula of the form $\fal x\,\alp(x)$ is a \textit{universal formula}.
In a proof-search of sequents in $\fal^{(++,-)}\exi^{(+)}$ [sequents in $\fal^{(++)}\exi^{(-)}$],
 only sequents in $\fal^{(++,-)}\exi^{(+)}$ [$\fal^{(++)}\exi^{(-)}$]
are produced, resp.

For a sequent $S$,
let $VC(S)$ denote the set of free variables and individual constants occurring in a formula in $S$.
If there is no such free variable nor individual constant in $S$,
put $VC(S)=\{c\}$ for an individual constant $c$.

A sequent $S=(\Gam\Rarw\Del)$ is \textit{saturated} iff 
it enjoys the five conditions in Definition \ref{df:saturated}, and it enjoys the following:
\benu
\item
if $(\exi y\,\bet(y))\in\Gam$, then $(\bet(b))\in\Gam^{\times}$ for a free variable $b\in FV$.
\item
if $(\exi y\,\bet(y))\in\Del$, then $(\bet(a))\in\Del^{\times}$ for every $a\in VC(S)$.
\item
if $(\fal x\,\alp(x))\in\Gam$, then $(\alp(a))\in\Gam^{\times}$ for every $a\in VC(S)$.
\eenu
A saturated sequent $\Gam\Rarw\Del$ is \textit{fully analyzed} 
if it is not an axiom, 
any marked formula $\alp^{\circ}$ in $\Gam\cup\Del$ is one of atomic formulas $R(t_{1},\ldots,t_{n})^{\circ}$,
$\bot^{\circ}$ ($\bot^{\circ}\not\in\Gam$).

Inversion steps for $\land,\lor,\supset$ are defined as in Definition \ref{df:inversionIp}.
Inversion steps for quantifiers are defined as follows.
\benu
\item
If an antecedent $\Gam$ of a sequent $S=(\Gam\Rarw\Del)$ contains a universal formula $(\fal y\,\bet(y))^{\circ}$ with a circle,
then erase the circle and add  starred instances for $VC(S)$:
for $\Gam=\{(\fal y\,\bet(y))^{\circ}\}\cup\Gam_{1}$ with $(\fal y\,\bet(y))^{\circ}\not\in\Gam_{1}$,
\[
\infer{(\fal y\,\bet(y))^{\circ},\Gam_{1}\Rarw\Del}
{
\{\bet(a)^{\st}:a\in VC(S)\}\cup\Gam_{1}\cup\{\fal y\,\bet(y)\}\Rarw\Del
}
\]
where $\{\bet(a)^{\st}:a\in VC(S)\}=\{\bet(a)^{\circ}: a\in VC(S), \bet(a)\not\in\Gam_{1}\}$.

\item
If a succedent $\Del$ of a sequent $S=(\Gam\Rarw\Del)$ contains a universal formula $(\fal y\,\bet(y))^{\circ}$ with a circle,
then erase the circle:
for $\Del=\Del_{1}\cup\{(\fal y\,\bet(y))^{\circ}\}$ with $(\fal y\,\bet(y))^{\circ}\not\in\Del_{1}$,
\[
\infer{\Gam\Rarw\Del_{1},(\fal y\,\bet(y))^{\circ}}
{
\Gam\Rarw\Del_{1},\fal y\,\bet(y)
}
\]

\item
If a succedent $\Del$ of a sequent $S=(\Gam\Rarw\Del)$ contains an existential formula $(\exi y\,\bet(y))^{\circ}$ with a circle,
then erase the circle and add instances for $VC(S)$ together with starred conclusions:
for $\Del=\Del_{1}\cup\{(\exi y\,\bet(y))^{\circ}\}$ with $(\exi y\,\bet(y))^{\circ}\not\in\Del_{1}$,
\[
\infer{\Gam\Rarw\Del_{1},(\exi y\,\bet(y))^{\circ}}
{
\Gam\Rarw\{\exi y\,\bet(y)\}\cup\Del_{1}\cup\{\bet(a),\bet(a)^{{\tt c}\st}:a\in VC(S)\}
}
\]
where
$\{\bet(a),\bet(a)^{{\tt c}\st}:a\in VC(S)\}=\{\bet(a):a\in VC(S)\}\cup\{\bet(a)^{{\tt c}\circ}:a\in VC(S), \bet(a)^{{\tt c}}\not\in\Del_{1}\}$.

\item
Suppose that an antecedent $\Gam$ contains an existential formula $(\exi y\,\bet(y))^{\circ}$ with a circle,
and $\Gam$ contains an instance $\bet(a)$ of $\exi y\,\bet(y)$ for a free variable $a$.
Then erase the circle:
Let $\Gam=\{(\exi y\,\bet(y))^{\circ},\bet(a)\}\cup\Gam_{1}$ with $(\exi y\,\bet(y))^{\circ}\not\in\Gam_{1}$
\[
\infer{(\exi y\,\bet(y))^{\circ},\bet(a),\Gam_{1}\Rarw\Del}
{
\Gam_{1},\exi y\,\bet(y),\bet(a)\Rarw\Del
}
\]

\item
Suppose that an antecedent $\Gam$ contains an existential formula $(\exi y\,\bet(y))^{\circ}$ with a circle,
and $\Gam$ contains no instance $\bet(a)$ of $\exi y\,\bet(y)$ for free variables $a$.
Then erase the circle, add a starred instance with an eigenvariable and instances for the eigenvariable together with starred conclusions:
Let $\Gam_{\fal}$ denote the set of universal unmarked formulas in $\Gam$,
and $\Del_{\exi}$ the set of existential unmarked formulas in $\Del$.
Then $\Gam=\{(\exi y\,\bet(y))^{\circ}\}\cup\Gam_{1}\cup\Gam_{\fal}$ with $(\exi y\,\bet(y))^{\circ}\not\in\Gam_{1}$
and $\Del=\Del_{1}\cup\Del_{\exi}$
\[
\infer{(\exi y\,\bet(y))^{\circ},\Gam_{1},\Gam_{\fal}\Rarw\Del_{1},\Del_{\exi}}
{
\bet(b)^{\circ},\Gam_{1},\exi y\,\bet(y),\Gam_{\fal}(b)^{\st}\Rarw\Del_{1},\Del_{\exi}(b),\Del_{\exi}(b)^{{\tt c} \st}
}
\]
where the variable $b$ does not occur in the lower sequent,
$\Gam_{\fal}(b)^{\st}=\{\alp(b)^{\circ}:\fal x\,\alp(x)\in\Gam_{\fal},\alp(b)\not\in\Gam\}$,
$\Del_{\exi}(b)=\{\gam(b):\exi z\,\gam(z)\in\Del_{\exi}\}$
and
$\Del_{\exi}(b)^{{\tt c}\st}=\{\gam(b)^{{\tt c}\circ}: \exi z\,\gam(z)\in\Del_{\exi}, \gam(b)^{{\tt c}}\not\in\Del\}$.

This inversion step for existential antecedent formulas is a condensed one.
The lower sequent is derivable from the upper one using some $(L\fal)$ and some $(R\exi)$ followed by
an $(L\exi)$.
\eenu
Observe that if the lower sequent enjoys the the condition (\ref{eq:saturated}),
then so do the upper sequents in each inversion step.

As far as one of the inversions can be performed, continue it to construct a tree $Tr_{S}$.
As in Proposition \ref{prp:Trterminate} we see that the process terminates, and $Tr_{S}$ is a finite tree.
Its depth is bounded by the size $\#S$ of the sequent $S$.
Assume that $S$ is a saturated sequent.
Then each leaf in $Tr_{S}$ is seen to be a saturated sequent, which is either
an axiom or a fully analyzed sequent.

Let $S=(\Gam,\Gam_{\supset},\Gam_{\fal}\Rarw \Del,\Del_{\supset}^{r},\Del_{\fal})$ be a sequent in 
$\fal^{(++,-)}\exi^{(+)}$,
where 
$\Gam_{\supset}$ denotes the set of unmarked implicational formulas in the antecedent,
$\Gam_{\fal}$ the set of unmarked universal formulas in the antecedent,
$\Del_{\supset}^{r}=S_{\supset}^{r}$ defined in Definition \ref{df:Sr}, 
and $\Del_{\fal}$ the set of (unmarked) universal formulas in the succedent, resp.
$\Gam\cup\Del$ is the remainders.
The set $\Gam_{\fal}$ is absent when
$S$ is a sequent in $\fal^{(++)}\exi^{(-)}$.

Assume that $\Del_{\supset}^{r}\cup\Del_{\fal}\neq\emptyset$.
The branching rule $(\mbox{br})$ here is of the form:
\[
\infer[(\mbox{br})]{\Gam,\Gam_{\supset},\Gam_{\fal}\Rarw\Del,\Del_{\supset}^{r},\Del_{\fal}}
{
\{\alp^{{\tt p}\st}\cup\Gam\cup\Gam_{\supset}^{\circ},\Gam_{\fal}^{\circ}\Rarw\alp^{{\tt c}\circ}: \alp\in\Del_{\supset}^{r}\}
& 
\{\Gam,\Gam_{\supset}^{\circ},\Gam_{\fal}^{\circ}\Rarw \gam(a),\gam(a)^{{\tt c}\circ}:(\fal y\,\gam(y))\in\Del_{\fal}\}
}
\]
where $a$ is an eigenvariable of the $(\mbox{br})$ and does
not to occur in the sequent $\Gam,\Gam_{\supset},\Gam_{\fal}\Rarw \fal y\,\gam(y)$.
Implicational formulas and universal formulas in the antecedent are marked to be analyzed again.
Each formula $\alp\in\Del_{\supset}^{r}$ and each $(\fal y\,\gam(y))\in\Del_{\fal}$ is a \textit{major formula}
of the rule $(\mbox{br})$,
and each $\gam(a)$ is a \textit{minor formula}.

Let $S_{0}$ be a given saturated sequent in one of classes $\fal^{(++,-)}\exi^{(+)}$ and $\fal^{(++)}\exi^{(-)}$ .
As in the propositional case, Definition \ref{df:TRp},
a $(\land,\lor)$-tree of deductions $TR(S_{0})$ is constructed from deductions $d(\sig)=Tr_{S(\sig)}$
for $\sig\in TR(S_{0})$ with $g(\sig)=\land$, where 
 $S(\sig)$ denotes the sequent at the node $\sig$ in $TR(S_{0})$.

\bdf
{\rm
For a formula $\alp$, let $d_{\fal}^{++}(\alp)$ denote
the number of nesting of strictly positive universal quantifiers in $\alp$.
The number is defined recursively as follows.
$d_{\fal}^{++}(\alp)=0$ if $\alp$ is an atomic formula.
$d_{\fal}^{++}(\alp_{0}\lor\alp_{1})=d_{\fal}^{++}(\alp_{0}\land\alp_{1})=\max\{d_{\fal}^{++}(\alp_{i}):i\in\{0,1\}\}$.
$d_{\fal}^{++}(\alp\supset\bet)=d_{\fal}^{++}(\bet)$.
$d_{\fal}^{++}(\exi x\,\alp)=d_{\fal}^{++}(\alp)$.
$d_{\fal}^{++}(\fal x\,\alp)=1+d_{\fal}^{++}(\alp)$.

Also let $d_{\fal}^{++}(S_{0}):=\max\{d_{\fal}^{++}(\alp):\alp\mbox{ is a succedent formula in } S_{0}\}$
for sequents $S_{0}$.
}
\edf

\subsection{The fragment $\fal^{(++)}\exi^{(-)}$}

\blem\label{prp:termination++}
The whole process generating the tree $TR(S_{0})$ for $\fal^{(++)}\exi^{(-)}$-sequents terminates, and
the number of $\lor$-gates along any branch in the tree $TR(S_{0})$ is at most 
$n_{\supset}^{+}(S_{0})+d_{\fal}^{++}(S_{0})$,
where
$n_{\supset}^{+}(S_{0})=\#Sbfml_{\supset}^{+}(S_{0})$ defined before Lemma \ref{prp:termination}.
\elem
\bprf

We say that a formula $\alp$ is an \textit{instance} of a formula $\bet$ if
$\alp$ is obtained from $\bet$ by replacing some variables
by variables and individual constants.

Each major formula of a branching rule in the tree $TR(S_{0})$ is an instance
of a positive subformula in $S_{0}$, which is either implicational or universal.

Consider a branching rule in the tree $TR(S_{0})$ one of whose major formulas
is a universal formula $\fal y\,\gam(y)$.
\[
\infer[(\mbox{br})]{\Gam,\Gam_{\supset}\Rarw\Del,\Del_{\supset}^{r},\Del_{\fal}}
{
\cdots
& 
\infer*{\Gam,\Gam_{\supset}^{\circ}\Rarw \gam(a),\gam(a)^{{\tt c}\circ}}
 {
  \infer[(\mbox{br})\, I]{S}{}
  }
&
\cdots
}
\]
where $(\fal y\,\gam(y))\in\Del_{\fal}$ and $a$ is the eigenvariable not occurring in 
$\Gam,\Gam_{\supset}\Rarw\fal y\,\gam(y)$.

Suppose that there occurs a branching rule $I$ above an upper sequent 
$S_{\gam}=(\Gam,\Gam_{\supset}^{\circ}\Rarw \gam(a),\gam(a)^{{\tt c}\circ})$
 one of whose major formulas
is a universal formula $\fal z\,\del(z)$.
Then the formula $\fal z\,\del(z)$ is an instance of a subformula of $\gam(a)^{{\tt c}}$.

Therefore each universal subformula in $S_{0}$ is analyzed at most once along any branch in $TR(S_{0})$,
and
there are at most $d_{\fal}^{++}(S_{0})$ applications of the branching rules along any branch
whose minor formula is an instance of the matrix of a universal major formula.
This means that there occurs at most one instance of the matrix $\gam(y)$ 
of a universal subformula $\fal y\,\gam(y)$ in $S_{0}$ along any branch in $TR(S_{0})$.

Next consider a negative existential formula $\exi y\,\bet(y)$.
Such a formula is analyzed in an inversion step
\[
\infer{(\exi y\,\bet(y))^{\circ},\Gam_{1},\Gam_{\fal}\Rarw\Del_{1},\Del_{\exi}}
{
\bet(b)^{\circ},\Gam_{1},\exi y\,\bet(y),\Gam_{\fal}(b)^{\st}\Rarw\Del_{1},\Del_{\exi}(b),\Del_{\exi}(b)^{{\tt c} \st}
}
\]
Once it is analyzed, it will not be analyzed again, since the formula $\bet(b)$ is in the antecedents
in sequents occurring above the upper sequent 
\\
$\bet(b)^{\circ},\Gam_{1},\exi y\,\bet(y),\Gam_{\fal}(b)^{\st}\Rarw\Del_{1},\Del_{\exi}(b),\Del_{\exi}(b)^{{\tt c} \st}$.
This means that there occurs at most one instance of the matrix $\bet(y)$ 
of an existential subformula $\exi y\,\bet(y)$ in $S_{0}$ along any branch in $TR(S_{0})$.

Since a new variable is introduced only when either a universal formula in succedents or an existential formula in
antecedents is analyzed,
each bound variable is replaced by at most one free variable along any branch.
Hence for each formula $\alp$ occurring in $S_{0}$,
only one instance of $\alp$ occurs in the whole tree $TR(S_{0})$ of deductions.
Therefore for each positive implicational subformula $\alp$ in $S_{0}$,
along any branch 
there is at most one instance $\alp^{\prime}$ of $\alp$ which is one of major formulas of a branching rule.
Furthermore as in Lemma \ref{prp:termination} we see that 
there is at most one branching rule along any branch one of whose major formula is the instance $\alp^{\prime}$.

Therefore the number of $\lor$-gates along any branch in the tree $TR(S_{0})$ is at most 
$n_{\supset}^{+}(S_{0})+d_{\fal}^{++}(S_{0})$.

Let $De(S_{0})$ denote the finite deduction with the branching rule $(\mbox{br})$, obtained from the finite tree
 $TR(S_{0})$ by carrying out intermediate deductions $Tr_{S(\sig)}$ for 
$\sig\in TR(S_{0})$ with $g(\sig)=\land$.
Since the size $\# S$ of sequents $S$ in $De(S_{0})$ is bounded by the square $(\#S_{0})^{2}$ of the size
$\#S_{0}$ of the given sequent $S_{0}$,
the depth of each deduction $Tr_{S}$ is also bounded by $(\#S_{0})^{2}$.
Hence the depth of the deduction $De(S_{0})$ is bounded by $(\#S_{0})^{3}$, and 
its width is bounded by $\#S_{0}$.

Therefore the whole process generating the tree $De(S_{0})$ terminates.
\eprf
\\

As for the propositional case, let us compute the value $v(\sig)$ of the node $\sig$ in the $(\land,\lor)$-tree $TR(S_{0})$.
If the value $v(\emptyset)$ is $1$, then $S_{0}$ is derivable.
Otherwise shrink the tree $TR(S_{0})$ to a tree $T\subset TR(S_{0})$ in a bottom-up manner.
Simultaneously a structure $I_{T}(\sig)$ with a universe $D_{T}(\sig)$ is assigned.
Suppose $v(\sig)=0$ for a node $\sig\in T$ with $g(\sig)=\land$.
Pick a son $\sig*(i)$ such that $v(\sig*(i))=0$, identify the node $\sig*(i)$ with $\sig$, and
choose the underivable sequent $S(\sig*(i))$, which is a non-axiom leaf
in the deduction $Tr_{S(\sig)}$.
Let $\Gam(\sig)\Rarw\Del(\sig)$ be the leaf sequent $S(\sig*(i))$ chosen from $Tr_{S(\sig)}$.
Then let $D_{T}(\sig)=VC(S(\sig*(i))$.
$R(a_{1},\ldots,a_{n})$ is defined to be true in the structure $I_{T}(\sig)$ iff $R(a_{1},\ldots,a_{n})\in \Gam(\sig)^{\times}$
for $a_{1},\ldots,a_{n}\in D_{T}(\sig)$ and $n$-ary predicate symbol $R$.
If $g(\sig*(i))=0$, then $\sig$ will be a leaf in $T$.
Otherwise $g(\sig*(i))=\lor$, and keep its sons in a shrunken tree, i.e., $\sig*(i,j)\in T$.

Thus we have defined a Kripke model $\langle T,\subset_{e},D_{T},I_{T}\rangle$.

For each $\sig\in T$, $\Gam(\sig)\Rarw\Del(\sig)$ denotes the leaf sequent
$S(\sig*(i))$ chosen from $Tr_{S(\sig)}$.

From the construction we see readily the followings for any $\sig,\tau\in T$.
\benu
\item
$\sig\subset_{e}\tau\Rarw \Gam(\sig)^{\times}\subset\Gam(\tau)^{\times}$.
\item
$\Gam(\sig)\Rarw\Del(\sig)$ is saturated.
\item
 if an implicational formula $\alp\in \Del(\sig)$, then there exists an extension $\tau\in T$ of $\sig$ such that
 $\alp^{{\tt p}}\subset \Gam(\tau)^{\times}$ and $\alp^{{\tt c}}\in\Del(\tau)^{\times}$.

\item
if a universal formula $\fal x\,\alp(x)\in\Del(\sig)$, then there exists an extension $\tau\in T$ of $\sig$ such that
$\alp(a)\in\Del(\tau)^{\times}$ for a variable $a$.
\item
$\Gam(\sig)^{\times}\cap\Del(\sig)^{\times}$ has no common atom, and $\bot\not\in\Gam(\sig)^{\times}$.
\eenu

\bprp\label{prp:adequacy++,-}
If $\alp\in\Gam(\sig)^{\times}$ $[\alp\in\Del(\sig)^{\times}]$, then $\sig\models\alp$ $[\sig\not\models\alp]$, resp. in the
 Kripke model $\langle T,\subset_{e},D_{T},I_{T}\rangle$.
Hence $\sig\models\bigwedge\Gam(\sig)$ and $\sig\not\models\bigvee\Del(\sig)$.
\eprp
\bprf
By simultaneous induction on $\alp$ using the above facts.
\eprf
\\

We conclude the followings from Proposition \ref{prp:adequacy++,-} and the proof of Lemma \ref{prp:termination++}.

\bth\label{th:intcomplete++,-}{\rm (Sch\"utte's dichotomy)}\\
For every saturated sequent $S_{0}$ in $\fal^{(++)}\exi^{(-)}$,
$v(\emptyset)=1$ iff ${\sf LJm}^{{\tt c}}\vdash S_{0}$.
\end{theorem}

\bcor\label{cor:depthIp++,-}
Each intuitionistically valid sequent $S_{0}$ in the class $\fal^{(++)}\exi^{(-)}$
has a derivation $\mathcal{D}$ in ${\sf LJm}^{{\tt c}}$ such that
$\mathcal{D}$ is a binary tree and
the depth of the tree is bounded by a cubic polynomial in the size $\#S_{0}$ of the sequent $S_{0}$.
\ecor

\bcor\label{cor:PSPACE++-}
The decision problem of the intuitionistic validity for formulas in the class $\fal^{(++)}\exi^{(-)}$ is 
solvable in PSPACE.
\ecor

\subsection{The fragment $\fal^{(++,-)}\exi^{(+)}$}

\bdf
{\rm
For a formula $\alp$, $q(\alp)$ denote the maximal number of bound variables occurring in a positive existential subformula,
and of bound variables occurring in a negative universal subformula in $\alp$.
It is defined recursively as follows.

$q^{\pm}(\alp)=0$ if $\alp$ is an atomic formula.
$q^{\pm}(\alp_{0}\lor\alp_{1})=q^{\pm}(\alp_{0}\land\alp_{1})=\max\{q^{\pm}(\alp_{i}):i=0,1\}$.
$q^{\pm}(\alp\supset\bet)=\max\{q^{\mp}(\alp),q^{\pm}(\bet)\}$.
$q^{+}(\fal x\,\alp)=q^{+}(\alp)$.
$q^{-}(\exi x\,\alp)=q^{-}(\alp)$.
\[
q^{-}(\fal x\,\alp)=
\left\{
\begin{array}{ll}
1+q^{-}(\alp) & \mbox{{\rm if }} x \mbox{ {\rm occurs in }} \alp
\\
q^{-}(\alp) & \mbox{{\rm otherwise}}
\end{array}
\right.
\]
\[
q^{+}(\exi x\,\alp)=\left\{
\begin{array}{ll}
1+q^{+}(\alp) & \mbox{{\rm if }} x \mbox{ {\rm occurs in }} \alp
\\
q^{+}(\alp) & \mbox{{\rm otherwise}}
\end{array}
\right.
\]

Then let $q(\alp)=q^{+}(\alp)$.

Also let $q(S):=\max(\{q^{-}(\alp):\alp\mbox{ is an antecedent formula in } S\}\cup\{q^{+}(\alp):\alp\mbox{ is a succedent formula in } S\})$
for sequents $S$.

}
\edf

\blem\label{prp:termination++,-}
The whole process generating the tree $TR(S_{0})$ for $\fal^{(++,-)}\exi^{(+)}$-sequents terminates, and
the number of $\lor$-gates along any branch in the tree $TR(S_{0})$ is at most 
$d_{\fal}+n\cdot d^{q}$,
where $n=\#Sbfml_{\supset}^{+}(S_{0})$, $d_{\fal}=d_{\fal}^{++}(S_{0})$, $d=d_{\fal}+\# VC(S_{0})$ and 
$q=q(S_{0})$.
\elem
\bprf
As in Lemma \ref{prp:termination++} we see that
each universal subformula in $S_{0}$ is analyzed at most once along any branch in $TR(S_{0})$,
and
there are at most $d_{\fal}$ applications of the branching rules along any branch
one of whose minor formula is an instance of the matrix of a universal major formula.
This means that there occurs at most one instance of the matrix $\gam(y)$ 
of a universal subformula $\fal y\,\gam(y)$ in $S_{0}$ along any branch in $TR(S_{0})$,
and the number of terms (free variables and individual constants) occurring in $De(S_{0})$ is bounded by $d$.

Hence for each formula $\alp$ occurring in $S_{0}$,
the number of instances of $\alp$ occurring in the whole tree $De(S_{0})$ is bounded by $d^{q}$.
Therefore for each positive implicational subformula $\alp$ in $S_{0}$,
along any branch 
there are at most $d^{q}$ instances $\alp^{\prime}$ of $\alp$ which is one of major formulas of a branching rule.
Furthermore as in Lemma \ref{prp:termination} we see that 
there is at most one branching rule along any branch one of whose major formula is the instance $\alp^{\prime}$.

Therefore the number of $\lor$-gates along any branch in the tree $TR(S_{0})$ is at most 
$d_{\fal}+n\cdot d^{q}$, and the whole process generating the tree $TR(S_{0})$ terminates.
\eprf
\\

As for propositional case and the the fragment $\fal^{(++)}\exi^{(-)}$ we see the followings.

\bth\label{th:intcomplete++,-}{\rm (Sch\"utte's dichotomy)}\\
For every saturated sequent $S_{0}$ in $\fal^{(++,-)}\exi^{(+)}$,
$v(\emptyset)=1$ iff ${\sf LJm}^{{\tt c}}\vdash S_{0}$.
\end{theorem}

\bcor\label{cor:depthIp++,-,+}
There exists a constant $c>0$ for which the following hold.

Each intuitionistically valid sequent $S_{0}$ in the class $\fal^{(++,-)}\exi^{(+)}$
has a derivation $\mathcal{D}$ in ${\sf LJm}^{{\tt c}}$ such that
$\mathcal{D}$ is a binary tree and
the depth of the tree is bounded by $cn^{n}$
for the size $n=\#(S_{0})$ of the sequent $S_{0}$.
Also the size of the sequents occurring in $\mathcal{D}$ is bounded by $cn^{n}$.
\ecor

\bcor\label{cor:PSPACE++-+}
The decision problem of the intuitionistic validity for formulas in the class $\fal^{(++,-)}\exi^{(+)}$ is 
solvable in exponential space.
\ecor

\section{Proof search in {\sf LJm}}\label{sec:Iqfull}

In this section we consider a proof search procedure for the full intuitionistic predicate logic {\sf Iq} 
over a finite language possibly with function symbols.
Here a search tree may be infinite.

{\sf LJm} denotes the standard sequent calculus {\sf m-G3i} for the intuitionistic predicate logic {\sf Iq} in \cite{TroSch}.
For example the right rule for existential formula is of the form.
\[
\infer[(R\exi)]{\Del\Rarw\Gam}
{
\Del\Rarw\Gam,\alp(t)
}
\]
where $\exi x\,\alp(x)\in\Gam$ and the minor formula $\alp(t)$ need not to be accompanied with its conclusion $\alp(t)^{{\tt c}}$.

For sequents $S$, the search tree $Tr_{S}$ is in general infinite due to the presence of \textit{universal formulas} $\fal x\,\alp(x)$
in antecedents and \textit{existential formulas} $\exi y\,\bet(y)$ in succedents.
A formula is \textit{non-invertible} if it is either an implicational formula or a universal formula.
It is desirable for us that each stage in constructing the tree of deductions is executed in a finite number of steps.
In order to do so, each stage tests only a finite number of free variables for universal formulas in antecedents and for existential formulas in succedents.
Let $\{t_{i}\}_{i}$ be an enumeration of all terms.
$Tm(A)$ denotes the set of all terms over a set $A\subset FV$ of free variables, and
\[
Tm(A)\restrict n:=\{t_{i}\in Tm(A): i<n\}
.\]

Let $n<\ome$ and $A$ a set of free variables.
A sequent $\Gam\Rarw\Del$ is $(n,A)$-\textit{saturated} iff 
it is saturated with respect to propositional connectives $\lor,\land$ as in Definition \ref{df:saturated} ,
and it enjoys the following four conditions:
\benu
\item if $\alp\supset\bet\in\Gam$, then $\alp\in\Gam^{\times}$ and $\bet\in\Del^{\times}$.
\item
if $(\exi y\,\bet(y))\in\Gam$, then $(\bet(a))\in\Gam^{\times}$ for a free variable $a\in FV$.
\item
if $(\exi y\,\bet(y))\in\Del$, then $(\bet(t))\in\Del^{\times}$ for every  $t\in Tm(A)\restrict n$.
\item
if $(\fal x\,\alp(x))\in\Gam$, then $(\alp(t))\in\Gam^{\times}$ for every  $t\in Tm(A)\restrict n$.
\eenu
An $(n,A)$-saturated sequent $\Gam\Rarw\Del$ is $(n,A)$-\textit{analyzed} 
if it is not an axiom, 
any marked formula $\alp^{\circ}$ in $\Gam\cup\Del$ is one of atomic formulas $R(t_{1},\ldots,t_{n})^{\circ}$,
$\bot^{\circ}$ ($\bot^{\circ}\not\in\Gam$),
or non-invertible formulas in $\Del$.
An $(n,A)$-analyzed sequent is \textit{fully analyzed} if
there is no existential succedent formula [no universal antecedent formula], resp.
Note that each fully analyzed sequent is $0$-analyzed.

A deduction $Tr_{S}^{(n,A)}$ is constructed in a finite number of steps 
as for propositional case by leaving any non-invertible succedent formulas
and applying $(L\exi), (R\exi), (L\fal)$ up to the $n$-th terms in $Tm(A)$
Put the given sequent $S$ at the root of the tree $Tr_{S}^{(n,A)}$.
The inversion steps for quantifiers are as follows.
\[
\infer[(\exi\Rarw)]{\exi x\,\alp(x)^{\circ},\Gam\Rarw\Del}
{\alp(a)^{\circ},\Gam,\exi x\,\alp(x)\Rarw\Del}
\]
where the eigenvariable $a$ 
does not occur in the lower sequent nor  nor in the finite set 
$A$ of free variables.
Moreover $a$ is chosen so that the condition (\ref{eq:freshvar}) in the next subsubsection \ref{subsec:extension} is met.
\[
\infer[(\Rarw\exi)]{\Del\Rarw\Gam,\exi x\,\alp(x)^{\circ}}
{
\Del\Rarw\exi x\,\alp(x),\Gam,\{\alp(t)^{\circ}\}_{t\in Tm(A)\restrict n}
}
\]
\[
\infer[(\fal\Rarw)]{\fal x\,\alp(x)^{\circ},\Gam\Rarw\Del}
{
\{\alp(t)^{\circ}\}_{t\in Tm(A)\restrict n},\Gam,\fal x\,\alp(x)\Rarw\Del}
\]
All of terms in the finite set $Tm(A)\restrict n$ are tested 
for existential formulas in succedent and for universal formulas in antecedent.

Each leaf in $Tr_{S}^{(n,A)}$ is $n$-saturated,which is either an axiom or an $(n,A)$-analyzed sequent
if $S$ is $(n,A)$-saturated.

\subsection{Extensions for non-invertible succedent formulas and postponed instantiations}\label{subsec:extension}

In a $\lor${\bf -stage} of our proof search for the predicate logic we examine all possibilities 
with succedent non-invertible formulas
by introducing a branching rule $(\mbox{br})$ as for the propositional case.
Consider an $n$-saturated sequent 
$$
\Gam,\Gam_{\supset},\Gam_{\fal}\Rarw\Del_{\exi},\Del,\Del_{\supset}^{\circ},\Del_{\fal}^{\circ}
$$
where 
$\Gam_{\supset}$ denotes the set of unmarked implicational formulas, and
$\Gam_{\fal}$ the set of unmarked universal formulas in the antecedent, resp.
$\Del_{\exi}$ denotes the set of unmarked existential formulas,
$\Del_{\supset}^{\circ}$ the set of marked implicational formulas, and
and $\Del_{\fal}^{\circ}$ the set of marked universal formulas in the succedent, resp.
$\Gam\cup\Del$ is the remainders.
Each marked formula in $\Gam\cup\Del$ is an atomic formula $R(t_{1},\ldots,t_{n})^{\circ}$ or $\bot^{\circ}$
with $\bot^{\circ}\not\in\Gam$.
Each unmarked formula in $\Gam$ is one of a disjunctive formula, a conjunctive formula
and an existential formula.
Each unmarked formula in $\Del$ is either disjunctive or conjunctive.

Assume that $\Gam_{\fal}\cup\Del_{\exi}\cup\Del_{\supset}^{\circ}\cup\Del_{\fal}^{\circ}\neq\emptyset$.

The sequent follows from several sequents.
Let us depict the several possibilities as an `inference rule' as follows.
{\scriptsize
\[
\infer[(\mbox{br})]{\Gam,\Gam_{\supset},\Gam_{\fal}\Rarw\Del_{\exi},\Del,\Del_{\supset}^{\circ},\Del_{\fal}^{\circ}}
{
\Gam,\Gam_{\supset},\Gam_{\fal}^{\circ}\Rarw\Del_{\exi}^{\circ},\Del,\Del_{\supset}^{\circ},\Del_{\fal}^{\circ}
&
\{\gam^{\circ},\Gam,\Gam_{\supset}^{\circ},\Gam_{\fal}^{\circ}\Rarw\del^{\circ}: (\gam\supset\del)^{\circ}\in\Del_{\supset}^{\circ}\}
& 
\{\Gam,\Gam_{\supset}^{\circ},\Gam_{\fal}^{\circ}\Rarw\gam(a)^{\circ}:(\fal y\,\gam(y))\in\Del_{\fal}^{\circ}\}
}
\]
}
where 
\benu

\item
for sets of unmarked formulas $\Gam$, $\Gam^{\circ}=\{\gam^{\circ}:\gam\in\Gam\}$.

The sequent $\Gam,\Gam_{\supset},\Gam_{\fal}^{\circ}\Rarw\Del_{\exi}^{\circ},\Del,\Del_{\supset}^{\circ},\Del_{\fal}^{\circ}$
is absent when $\Gam_{\fal}\cup\Del_{\exi}=\emptyset$.

\item
$a$ is an eigenvariable distinct each other for universal formulas $(\fal y\,\gam(y))\in\Del_{\fal}$
and such that the condition (\ref{eq:freshvar}) below is met.
\eenu
If one of upper sequents of $(\mbox{br})$ is derivable, 
then so is the lower sequent possibly using a non-invertible inference rule $(R\supset)$
or $(R\fal)$.
Each upper sequent $\gam^{\circ},\Gam,\Gam_{\supset}^{\circ},\Gam_{\fal}^{\circ}\Rarw\del^{\circ}$ with $(\gam\supset\del)^{\circ}\in\Del_{\supset}^{\circ}$
and each $\Gam,\Gam_{\supset}^{\circ},\Gam_{\fal}^{\circ}\Rarw\gam(a)^{\circ}$ with $(\fal y\,\gam(y))\in\Del_{\fal}^{\circ}$ is said to be a
\textit{non-invertible upper sequent} of the inference rule $(\mbox{br})$.
While the leftmost upper sequent $\Gam,\Gam_{\supset},\Gam_{\fal}^{\circ}\Rarw\Del_{\exi}^{\circ},\Del,\Del_{\supset}^{\circ},\Del_{\fal}^{\circ}$
is the \textit{continued} sequent.

Here is the condition on eigenvariables.
\beqnarr
&&
\mbox{
Each eigenvariable is distinct each other, and}
\nonumber
\\
&&
\mbox{ occurs only above the inference rule or}
\label{eq:freshvar}
\\
&&
\mbox{in the right part of the inference where the variable is introduced.}
\nonumber
\eeqnarr
This means that if $\sig$ is the node of the upper sequent of the inference rule where
an eigenvariable $a$ is introduced,
and $a$ occurs in the sequent at a node $\tau$,
then either $\sig\subset_{e}\tau$(, i.e., $\tau$ is above $\sig$) or
$\rho*(i)\subset_{e}\sig$ and $\rho*(j)\subset_{e}\tau$ for some $\rho$ and $i<j$
(, i.e., $\tau$ is right to $\sig$).

As in the propositional case, Definition \ref{df:TRp}, 
let us construct a tree $TR(S_{0})\subset{}^{<\ome}\ome$ for a given sequent $S_{0}$.
The tree $TR(S_{0})$ 
is constructed in $\ome$-steps.
$TR(S_{0})_{n}$ denotes the piece of $TR(S_{0})$ in the $n$th step
such that for each $\sig\in TR(S_{0})_{n}$,
the \textit{length} $lh(\sig)\leq n$.
The labeling function
$(S(\sig),d(\sig),g(\sig))$ for $\sig\in \bigcup_{n\in\ome}TR(S_{0})_{n}$ is defined simultaneously in the construction of $TR(S_{0})$.

$TR(S_{0})$ is defined to be 
the union $\bigcup_{n\in\ome}TR(S_{0})_{n}$.

Let $S(\tau)$ be the sequent at the node $\tau$ in $TR(S_{0})$, 
and $FV(\tau)$ the set of free variables occurring in the sequent $S(\tau)$.
A finite set $FV_{\subset_{e}}(\sig)$ of free variables is assigned to sequences $\sig\in TR(S_{0})$
as follows. The set is finite since the tree $TR(S_{0})_{n}\subset{}^{<n+1}\ome$ is finitely branching.
\[
FV_{\subset_{e}}(\sig)=\bigcup\{FV(\tau): \exi\rho(\rho\subset_{e}\sig \spand \rho\subset_{e}^{0}\tau\in TR(S_{0})_{lh(\sig)})\}
\] 
where
\[
\sig\subset^{0}_{e}\tau :\Lrarw \sig\subset_{e}\tau \spand 
\fal\rho[\sig\subset_{e}\rho\subsetneq_{e}\tau \spand g(\rho)=\lor
\Rarw \tau(lh(\rho))=0]
\]
with the $i$-th \textit{component} $\tau(i)$ of sequences $\tau$ for $i<lh(\tau)$.
$\sig\subset^{0}_{e}\tau$ means that $\tau$ continues to substitute terms for $\fal y$ in antecedents and $\exi x$ in succedent
for any $\lor${\bf -stage} after $\sig$.

Let us denote 
$$
Tr_{\sig}:=Tr_{S(\sig)}^{(lh(\sig),FV_{\subset_{e}}(\sig))}.
$$

$De(S_{0})$ denotes the whole tree of deductions obtained from $TR(S_{0})$ by fulfilling intermediate deductions,
and is constructed recursively.
Each $\land${\bf -stage} analyzes the current leaves parallel as in the propositional case.
After the $\land${\bf -stage}, we extend the tree by non-invertible 
$(\mbox{br})$ inference rules in $\lor${\bf -stage}.
In each moment $De(S_{0})$ is constructed so that 
the condition (\ref{eq:freshvar}) on eigenvariables is met.

\bdf\label{df:TRq}
{\rm Given a sequent $S_{0}=(\Gam_{0}\Rarw\Del_{0})$, let us define trees 
$TR(S_{0})_{n}\subset{}^{<n+1}\ome$, 
 and a labeling function $(S(\sig), d(\sig),g(\sig))$ for $\sig\in TR(S_{0})_{n}$,
where $S(\sig)$ is a sequent, $g(\sig)\in\{\lor,\land,0,1\}$ is a gate and $d(\sig)$ is a deduction possibly with the branching rule.

First the empty sequence $\emptyset\in TR(S_{0})_{0}=\{\emptyset\}$ and $S(\emptyset)=S_{0}$  where each formula in 
$S_{0}$ is marked.
Let $FV(\emptyset)$ be the set of free variables occurring in $S_{0}$ if the set is non-empty.
Otherwise $FV(\emptyset)=\{a_{0}\}$.
If $S_{0}$ is an axiom, then $g(\emptyset)=1$.
If $S_{0}$ is fully analyzed, then $g(\emptyset)=0$.
If $g(\emptyset)\in\{0,1\}$, then $d(\emptyset)$ is the deduction consisting solely of $S_{0}$.
Otherwise $g(\emptyset)=\land$ and the tree is extended according to the $\land${\bf -stage}.

Suppose that $TR(S_{0})_{n}$ has been constructed, and
there exists a leaf $\sig\in TR(S_{0})_{n}$ such that $g(\sig)\in\{\land,\lor\}$.
(Otherwise we are done, and $TR(S_{0})_{n+1}$ is not defined.)
If $n$ is even [odd], the tree is extended according to the $\land${\bf -stage} [the $\lor${\bf -stage}], resp.
\\

\noindent
$\land${\bf -stage}.

Consider each leaf $\sig\in TR(S_{0})_{n}$ with $g(\sig)=\land$. 
Extend the tree $De(S_{0})$ by putting the deduction $d(\sig)=Tr_{\sig}$ for each such $\sig$.
Let $\{S_{i}\}_{i<I}$ be an enumeration of all leaves in $d(\sig)$.
For each $i<I$, let $\sig*(i)\in TR(S_{0})_{n+1}$ with $S(\sig*(i))=S_{i}$.
If $S_{i}$ is an axiom, then $g(\sig*(i))=1$.
If $S_{i}$ is fully analyzed, then $g(\sig*(i))=0$.
Otherwise $S_{i}$ is not fully analyzed, but $(lh(\sig),FV_{\subset_{e}}(\sig))$-analyzed.
This means that either its antecedent contains a universal formula, or its succedent contains either an existential formula or a non-invertible formula.
Let $g(\sig*(i))=\lor$, and the tree is extended according to the $\lor${\bf -stage}.

$TR(S_{0})_{n+1}$ is defined to be the union of $TR(S_{0})_{n}$ and nodes $\sig*(i)$ for each leaf $\sig\in TR(S_{0})_{n}$
such that $g(\sig)=\land$ and $i<I$.
\\

\noindent
$\lor${\bf -stage}.

Consider each leaf $\sig\in TR(S_{0})_{n}$ with $g(\sig)=\lor$.
Extend the tree $De(S_{0})$ by the inference rule $(\mbox{br})$ parallel for each such $\sig$.

Let $S(\sig)=S$.
$S$ is a non-invertible sequent $\Gam,\Gam_{\supset},\Gam_{\fal}\Rarw\Del_{\exi},\Del,\Del_{\supset}^{\circ},\Del_{\fal}^{\circ}$
where cedents $\Gam,\Gam_{\supset},\Gam_{\fal},\Del_{\exi},\Del,\Del_{\supset}^{\circ}$ and $\Del_{\fal}^{\circ}$
are defined as in the beginning of this subsection, and
$\Gam_{\fal}\cup\Del_{\exi}\cup\Del_{\supset}^{\circ}\cup\Del_{\fal}^{\circ}\neq\emptyset$.
Let $S(\sig*(0))$ be the sequent
$\Gam,\Gam_{\supset},\Gam_{\fal}^{\circ}\Rarw\Del_{\exi}^{\circ},\Del,\Del_{\supset}^{\circ},\Del_{\fal}^{\circ}$.
Let $\Del_{\supset}^{\circ}=\{\bet_{j}^{\circ}\}_{0<j\leq J_{\supset}}$,
and
$\Del_{\fal}^{\circ}=\{\bet_{j}^{\circ}\}_{J_{\supset}<j\leq J_{\fal}}$.
Then $\sig*(j)\in TR(S_{0})_{n+1}$ for each $j$ with $0\leq j\leq J_{\fal}$.
For $j>0$ the sequent $S(\sig*(j))$ is defined by analyzing the $j$-th non-invertible formula $\bet_{j}$.
Namely if $\bet_{j}\equiv(\gam\supset\del)$, then
$S(\sig*(j))=(\gam^{\circ},\Gam,\Gam_{\supset}^{\circ},\Gam_{\fal}^{\circ}\Rarw\del^{\circ})$.
If $\bet_{j}\equiv(\fal y\,\gam(y))$, then
$S(\sig*(j))=(\Gam,\Gam_{\supset}^{\circ},\Gam_{\fal}^{\circ}\Rarw\gam(a)^{\circ})$
where the eigenvariables $a$ are fresh, i.e., do not occur in any $S(\tau)$ for $\tau\in TR(S_{0})_{n}$, and 
distinct each other for $(\fal y\,\gam(y))^{\circ}\in\Del_{\fal}$.
$d(\sig)$ denotes the deduction consisting of a $(\mbox{br})$ with its lower sequent $S=S(\sig)$.

Let $g(\sig*(j))=1$ if $S(\sig*(j))$ is an axiom, and $g(\sig*(j))=0$ if it is fully analyzed.
Otherwise let $g(\sig*(j))=\land$ and the tree is extended according to the $\land${\bf -stage}.

$TR(S_{0})_{n+1}$ is defined to be the union of $TR(S_{0})_{n}$ and nodes $\sig*(j)$ for each leaf $\sig\in TR(S_{0})_{n}$
such that $g(\sig)=\lor$ and $j\leq J_{\fal}$.

Finally let

\[
TR(S_{0})= \bigcup_{n\in\ome}TR(S_{0})_{n}
.\]

}
\edf

Let $\emptyset\neq T\subset TR(S_{0})$ be a subtree of $TR(S_{0})$ such that $g(\sig)\neq 1$ for any $\sig\in T$.
Let us construct a Kripke model
 $\langle T,\subset_{e},D_{T},I_{T}\rangle$
 as follows.

\bdf\label{df:Kripkepredicate}
\benu

\item
{\rm For} $S(\sig)=(\Gam(\sig)\Rarw\Del(\sig))$, {\rm let}
\beqnarrs
\Gam^{\infty}(\tau;T) & = & \bigcup\{\Gam(\rho)^{\times}:\tau\subset^{0}_{e}\rho\in T\}
\\
\Gam^{\infty}_{\subset_{e}}(\sig;T) & = & \bigcup\{\Gam^{\infty}(\tau;T): \tau\subset_{e}\sig\}
\\
\Del^{\infty}(\tau;T) & = & \bigcup\{\Del(\rho)^{\times}:\tau\subset^{0}_{e}\rho\in T\}
\eeqnarrs
{\rm Note that} $\Gam(\sig)^{\times}\subset\Gam^{\infty}(\sig;T)\subset\Gam^{\infty}_{\subset_{e}}(\sig;T) $ 
{\rm and} $\Del(\sig)^{\times}\subset\Del^{\infty}(\sig;T)$ {\rm for} $\sig\in T$.

\item
{\rm Let}
\[
D_{T}(\sig)=Tm(FV^{\infty}_{\subset_{e}}(\sig;T))
\]
{\rm where} 
\beqnarrs
FV^{\infty}(\tau;T) & = & \bigcup\{FV(\rho): \tau\subset^{0}_{e}\rho\in T\}
\\
FV^{\infty}_{\subset_{e}}(\sig;T) & = & \bigcup\{FV^{\infty}(\tau;T) : \tau\subset_{e}\sig\}
\eeqnarrs

\item
{\rm For} $\sig\in T$ {\rm and an $n$-ary predicate symbol} $R$ {\rm and function symbol} $f$,
{\rm let} 
\[
R^{\sig}=\{(t_{1},\ldots,t_{n}): t_{1},\ldots,t_{n}\in D_{T}(\sig)\spand R(t_{1},\ldots,t_{n})\in\Gam^{\infty}_{\subset_{e}}(\sig;T)\}
\]
{\rm and}
 $f^{\sig}(t_{1},\ldots,t_{n})=f(t_{1},\ldots,t_{n})$ {\rm for} $t_{1},\ldots,t_{n}\in D_{T}(\sig)$.
 
 \item
 $I_{T}(\sig)$ {\rm is a structure with its universe} $D_{T}(\sig)$ {\rm and relations} $R^{\sig}$ {\rm and functions} $f^{\sig}$.
\eenu
\edf

\bprp\label{prp:Kripkemodeldf}
Let $\emptyset\neq T\subset TR(S_{0})$ be a subtree of $TR(S_{0})$ such that $g(\sig)\neq 1$ for any $\sig\in T$.
Then
$\langle T,\subset_{e},D_{T},I_{T}\rangle$ is a Kripke model.
\eprp
\bprf
Let $\sig\subset_{e}\tau$ for $\sig,\tau\in T$.
Then $FV^{\infty}_{\subset_{e}}(\sig;T)\subset FV^{\infty}_{\subset_{e}}(\tau;T)$, and hence
$D_{T}(\sig)\subset D_{T}(\tau)$.
Moreover we have $\Gam^{\infty}_{\subset_{e}}(\sig;T)\subset\Gam^{\infty}_{\subset_{e}}(\tau;T)$,
and $R^{\sig}\subset R^{\tau}$.
\eprf
\\

A pair $\Gam\Rarw\Del$ of (possibly infinite) sets $\Gam,\Del$ of formulas is $A$\textit{-saturated} for a set $A$ of free variables
iff it is $(n,A)$-saturated for any $n$.
This means besides the saturation with respect to propositional connectives and existential formulas in $\Gam$ that 
\benu
\item
if $(\exi y\,\bet(y))\in\Del$, then $(\bet(t))\in\Del^{\times}$ for any $t\in Tm(A)$.
\item
if $(\fal x\,\alp(x))\in\Gam$, then $(\alp(t))\in\Gam^{\times}$ for any $t\in Tm(A)$.
\eenu

$A$-saturated pair $\Gam\Rarw\Del$ is
$A$\textit{-analyzed} 
if $\bot\not\in\Gam^{\times}$, $\Gam$ and $\Del$ has no common atomic formula.

\bprp\label{prp:countermodel}
Let $\emptyset\neq T\subset TR(S_{0})$ be a subtree of $TR(S_{0})$ such that $g(\sig)\neq 1$ for any $\sig\in T$.
Suppose that $T$ enjoys the following conditions for any $\sig\in T$.

\benu

\item
 $\Gam^{\infty}_{\subset_{e}}(\sig;T)\Rarw\Del^{\infty}(\sig;T)$ is 
$FV^{\infty}_{\subset_{e}}(\sig;T)$-analyzed.

\item
 \benu
 \item
 if $(\alp\supset\bet)\in \Del^{\infty}(\sig;T)$, then there exists an extension $\tau\in T$ of $\sig$ such that
 $\alp\in \Gam_{\subset_{e}}^{\infty}(\tau;T)$ and $\bet\in\Del^{\infty}(\tau;T)$.
 \item
  if $(\fal x\,\alp(x))\in \Del^{\infty}(\sig;T)$, then there exist an extension $\tau\in T$ of $\sig$ and 
  an $a\in FV_{\subset_{e}}^{\infty}(\tau;T)$
  such that
 $(\alp(a))\in\Del^{\infty}(\tau;T)$.
 
 \item
 $\Gam^{\infty}_{\subset_{e}}(\sig;T)\cap\Del^{\infty}(\sig;T)$ has no common atomic formula.
 \eenu
\eenu 

Let $\sig\in T$ 
and $\alp$ be a formula all of whose free variables are in the set $D_{T}(\sig)$.
In the
 Kripke model $\langle T,\subset_{e},D_{T},I_{T}\rangle$, 
 if $\alp\in\Gam^{\infty}_{\subset_{e}}(\sig;T)$, then $\sig\models\alp$,
 and if $\alp\in\Del^{\infty}(\sig;T)$, then $\sig\not\models\alp$.
 
Hence $\sig\models\bigwedge\Gam(\sig)$ and $\sig\not\models\bigvee\Del(\sig)$.
\eprp
\bprf
This is shown by simultaneous induction on $\alp$.

Consider the case when $\alp$ is an atomic formula $R(t_{1},\ldots,t_{n})$.
By the assumption we have $\alp\not\in\Gam^{\infty}_{\subset_{e}}(\sig;T)\cap\Del^{\infty}(\sig;T)$.
Hence if $\alp\in \Del^{\infty}(\sig;T)$, then $\alp\not\in\Gam^{\infty}_{\subset_{e}}(\sig;T)$, i.e.,
$\sig\not\models\alp$.
On the other side if $\alp\in\Gam^{\infty}_{\subset_{e}}(\sig;T)$
and $t_{1},\ldots,t_{n}\in D_{T}(\sig)$ by the assumption,
then $\sig\models\alp$.

Next consider the case when $\alp\equiv(\fal x\,\bet(x))$.
Suppose $\alp\in\Gam^{\infty}_{\subset_{e}}(\sig;T)$. 
For any extension $\tau$ of $\sig$ in $T$, i.e., $\sig\subset_{e}\tau\in T$,
$\alp\in\Gam^{\infty}_{\subset_{e}}(\sig;T)\subset\Gam^{\infty}_{\subset_{e}}(\tau;T)$.
$\bet(t)\in\Gam^{\infty}_{\subset_{e}}(\tau;T)$ for any $t\in D_{T}(\tau)$
by the supposition.
By IH $\tau\models\bet(t)$.
Hence $\sig\models\alp$.
Next suppose $\alp\in\Del^{\infty}(\sig;T)$.
Then by the supposition,
 for an extension $\tau\in T$ of $\sig$ and a free variable $a\in FV^{\infty}_{\subset_{e}}(\tau;T)$ , we have 
$\bet(a)\in\Del^{\infty}(\tau;T)$.
Hence $a\in D_{T}(\tau)$ as long as $a$ occurs in $\bet(a)$, and $\tau\not\models\bet(a)$ by IH.
Thus $\sig\not\models\alp$.

Other cases are seen easily.
\eprf
\\

The first condition (analyzed) and second one (existences of extensions) in Proposition \ref{prp:countermodel}
are easily enjoyed when $T$ has sufficiently many nodes, i.e., 
when nodes are prolonged in $T$ unlimitedly.
The third condition (no common atomic formula) is hard to satisfy.
Obviously the tree $TR(S_{0})$ of the deductions enjoys the first and second conditions.
But it may be the case that for a $\sig\in TR(S_{0})$ with $g(\sig)=\lor$,
$\Gam^{\infty}_{\subset_{e}}(\sig;TR(S_{0}))\cap\Del^{\infty}(\sig;TR(S_{0}))$ has a common atomic formula.

\bprp\label{prp:fanlyzedLJm}
Let $\emptyset\neq T\subset TR(S_{0})$ be a subtree of $TR(S_{0})$ such that $g(\sig)\neq 1$ for any $\sig\in T$.
Suppose that each $\land$-gate has a unique $\lor$-gate son in $T$.
Then
$\Gam^{\infty}(\sig;T)\cap\Del^{\infty}(\sig;T)$ has no common atomic formula for any $\sig\in T$.
\eprp
\bprf
Suppose that $\Gam^{\infty}(\sig;T)\cap\Del^{\infty}(\sig;T)$ has a common atomic formula $\alp$.
Let $\sig\subset_{e}^{0}\rho_{a},\rho_{s}\in T$
be such that $g(\rho_{a})\neq 1, g(\rho_{s})\neq 1$ and $\alp\in\Gam^{\times}(\rho_{a})\cap\Del^{\times}(\rho_{s})$.
We see that $\rho_{a}$ and $\rho_{s}$ are comparable in the order $\subset_{e}^{0}$,
since 
each $\land$-gate has a unique $\lor$-gate son in $T$,
and each $\lor$-gate has a unique son, i.e., the leftmost continued one in $\subset^{0}_{e}$.
Then $\alp\in\Gam^{\times}(\rho)\cap\Del^{\times}(\rho)$ for a common extension $\rho\in \{\rho_{a},\rho_{s}\}$,
and $S(\rho)$ is an axiom.
This means a contradiction $1=g(\rho)\neq 1$.
\eprf

\subsection{Completeness}\label{subsec:completeness}

In order to have the Sch\"utte's dichotomy,
we need to transform the tree $TR(S_{0})$ of deductions, cf.\,subsection \ref{subsec:transfer}.
However when we need only to show the completeness of {\sf LJm} with the cut rule $(cut)$, ${\sf LJm}+(cut)$,
one can extract a consistent tree $T$ from $TR(S_{0})$ by which $S_{0}$ is refuted
provided that $S_{0}$ is not derivable in ${\sf LJm}+(cut)$.

In this subsection we consider the tree $TR(S_{0})$.

For formulas $\alp$, $\alp^{\fal}$ denotes ambiguously formulas obtained from $\alp$ by binding some (possibly none)
free variables by universal quantifiers.

For formulas $\alp$ and $\bet$,
$\oplus(\alp,\bet)$ denotes the formula $\fal\vec{x}(\gam(\vec{x})\supset (\del(\vec{x})\lor\bet))$ if $\alp$ is a universally bound implicational formula 
$\fal\vec{x}(\gam(\vec{x})\supset \del(\vec{x}))$ for a list (possibly empty) list $\vec{x}$ of bound variables, where
$\vec{x}$ does not occur in $\bet$.
Otherwise $\oplus(\alp,\bet):\equiv (\alp\lor\bet)$.
For sequences $\vec{\alp}=\alp_{0},\alp_{1},\ldots,\alp_{n-1}$ of formulas $\alp_{i}$,
let
$\oplus(\vec{\alp},\bet):\equiv \oplus(\alp_{0},\oplus(\alp_{1},\cdots,\oplus(\alp_{n-1},\bet)\cdots))$, and
$\oplus^{\fal}(\vec{\alp},\bet)$ denotes formulas 
$\oplus^{\fal}(\alp_{0},\oplus^{\fal}(\alp_{1},\cdots,\oplus^{\fal}(\alp_{n-1},\bet)\cdots))$.

For the sequent $S(\sig)=(\Gam(\sig)\Rarw\Del(\sig))$ let
\[
\chi(\sig)  :\Lrarw  (\bigwedge \Gam(\sig)\supset  \bigvee\Del(\sig))
\]

Let $T\subset TR(S_{0})$ be a \textit{finite} subtree such that
for each $\sig\in T$, if $g(\sig)=\land$, then $\sig$ has a single son $\sig*(i)$ in $T$.
And if $g(\sig)=\lor$, then either $\sig$ is a leaf in $T$, or 
$\sig$ has all of sons $\sig*(j)$ in $T$, i.e., $\fal j[\sig*(j)\in T\lrarw \sig*(j)\in TR(S_{0})]$.
Moreover $g(\sig)\in\{0,\lor\}$ for any leaves $\sig$ in $T$.
Such a tree $T$ is called a \textit{selected} tree.
Following \cite{Kripke} let us introduce \textit{characteristic formulas} $\chi(\sig;T)$ for nodes $\sig$ of such a tree $T$
 recursively as follows.
Recall that $FV(\sig)$ denotes the set of free variables occurring in the sequent $S(\sig)$ for $\sig \in TR(S_{0})$.

For leaves $\sig$ in $T$,
\[
\chi(\sig;T)\equiv\chi(\sig)
.\]
Let $\sig\in T$ be an internal node with
$g(\sig)=\land$, and $\sig*(i)$ be the unique son of $\sig$ in $T$.
Also let $\chi(\sig*(i);T)\equiv \alp(\vec{a})$ where $\vec{a}$ is the set of eigenvariables of inference rules $(\exi\Rarw)$
occurring between the leaf $S(\sig*(i))$ and the root $\sig$ of the deduction $Tr_{\sig}$.
Then let
\[
\chi(\sig;T):\equiv \chi^{\fal}(\sig*(i);T) :\equiv  \fal \vec{x}\, \alp(\vec{x})
.\]
Let $\sig\in T$ be an internal node with
$g(\sig)=\lor$, and $\sig*(j)\,(j\leq J_{\fal})$ be all of sons of $\sig$ in $T$,
where $\Del_{\fal}^{\circ}=\{\bet_{j}\}_{J_{\supset}<j\leq J_{\fal}}$ for the set $\Del_{\fal}^{\circ}$ of marked universal formulas in
the succedent of $S(\sig)$.
For each $j$ with $J_{\supset}<j\leq J_{\fal}$
let $ \chi(\sig*(j);T)\equiv \alp(a)$ for the eigenvariable $a$ introduced at the $j$-th upper sequent $S(\sig*(j))$.
Then 
for $\chi^{\fal}(\sig*(j);T)\equiv \fal x\, \alp(x)$,
let
\[
\chi(\sig;T)\equiv \oplus(\chi(\sig*(0);T),\bigvee_{0<j\leq J_{\supset}}\chi(\sig*(j);T)\lor\bigvee_{J_{\supset}<j\leq J_{\fal}}\chi^{\fal}(\sig*(j);T))
.\]
Finally let
$\chi(T):\equiv\chi(\sig;T)$
for the root $\sig$ in $T$.

\bprp\label{prp:oplus}
\benu
\item\label{prp:oplus1}
${\sf LJm}+(cut)\vdash\alp\lor\bet \supset  \oplus(\alp,\bet)$.
${\sf LJm}+(cut)\vdash\oplus(\alp,\bot)\lrarw\alp\lrarw\oplus(\bot,\alp)$
and
${\sf LJm}+(cut)\vdash\oplus^{\fal}(\oplus^{\fal}(\alp,\bet),\gam)\lrarw\oplus^{\fal}(\alp,\bet\lor\gam)$.

\item\label{prp:oplus15}
Let $T\subset TR(S_{0})$ be a selected tree, and $\sig$ be a leaf in $T$.
Then there exist formulas $\vec{\alp}$ such that
${\sf LJm}+(cut)\vdash\chi(T)\lrarw\oplus^{\fal}(\vec{\alp},\chi(\sig))$.

\item\label{prp:oplus2}
Let $g(\sig)=\land$ and ${\sf LJm}+(cut)\not\vdash\oplus^{\fal}(\vec{\alp},\oplus(\chi^{\fal}(\sig),\bet))$ for some formulas $\vec{\alp},\bet$.
Then there exists a leaf $\sig*(i)$ in the deduction $Tr_{\sig}$ such that
${\sf LJm}+(cut)\not\vdash\oplus^{\fal}(\vec{\alp},\oplus(\chi^{\fal}(\sig*(i)),\bet))$.

\item\label{prp:oplus6}
Let $g(\sig)=\lor$ and ${\sf LJm}+(cut)\not\vdash\oplus^{\fal}(\vec{\alp},\oplus(\chi(\sig),\bet))$ for formulas $\vec{\alp},\bet$.
Then for each $j$ there exists a leaf $\sig*(j,i_{j})$ in the deduction $Tr_{\sig*(j)}$ such that
${\sf LJm}+(cut)\not\vdash\oplus(\vec{\alp},\oplus(\oplus(\chi(\sig*(0,i_{0})),\bigvee_{0<j\leq J_{\supset}}\chi(\sig*(j,i_{j}))\lor \bigvee_{J_{\supset}<j\leq J_{\fal}}\chi^{\fal}(\sig*(j,i_{j}))),\bet))$.

\item\label{prp:oplus7}
Let $T\subset TR(S_{0})$ be a selected tree.
Assume ${\sf LJm}+(cut)\not\vdash\chi(T)$.
Then for each leaf $\sig\in T$ and each $j$ there exists a leaf $\sig*(j,i_{j})$ in the deduction $Tr_{\sig*(j)}$ such that
${\sf LJm}+(cut)\not\vdash\chi(T^{\prime})$
for the tree $T^{\prime}$ obtained from $T$ by extending each leaf $\sig$ to $\sig*(j),\sig*(j,i_{j})$.
\eenu
\eprp
\bprf
\ref{prp:oplus}.\ref{prp:oplus15} is seen by induction on the size of the tree $T$
using Proposition \ref{prp:oplus}.\ref{prp:oplus1}.
\\

\noindent
\ref{prp:oplus}.\ref{prp:oplus2} is seen by inspection to inference rules in {\sf LJm} except non-invertible ones $(R\supset)$ and $(R\fal)$.
\\

\noindent
\ref{prp:oplus}.\ref{prp:oplus6}.
Assume $g(\sig)=\lor$ and ${\sf LJm}+(cut)\not\vdash\oplus(\vec{\alp},\oplus(\chi(\sig),\bet))$.
Then 
\[
{\sf LJm}+(cut)\not\vdash\oplus(\vec{\alp},\oplus(\oplus(\chi(\sig*(0)),\bigvee_{j>0}\chi^{\fal}(\sig*(j))),\bet))
\]
by the definition of the rule $(\mbox{br})$.
In other words,
\[
{\sf LJm}+(cut)\not\vdash\oplus(\vec{\alp},\oplus(\chi(\sig*(0)),\bigvee_{j>0}\chi^{\fal}(\sig*(j))\lor\bet))
.\]
By  Proposition \ref{prp:oplus}.\ref{prp:oplus2} pick an $i_{0}$ such that
\[
{\sf LJm}+(cut)\not\vdash\oplus(\vec{\alp},\oplus(\chi^{\fal}(\sig*(0,i_{0})),\bigvee_{j>0}\chi^{\fal}(\sig*(j))\lor\bet))
.\]
Hence
\[
{\sf LJm}+(cut)\not\vdash\oplus(\vec{\alp},\oplus(\oplus(\chi^{\fal}(\sig*(0,i_{0})),\bigvee_{j>1}\chi^{\fal}(\sig*(j)))\lor\bet,\chi(\sig*(1))))
.\]
Then again by Proposition \ref{prp:oplus}.\ref{prp:oplus2}
pick an $i_{1}$ such that
\[
{\sf LJm}+(cut)\not\vdash\oplus(\vec{\alp},\oplus(\oplus(\chi^{\fal}(\sig*(0,i_{0})),\bigvee_{j>1}\chi^{\fal}(\sig*(j)))\lor\bet,\chi^{\fal}(\sig*(1,i_{1}))))
.\]
In this way we can pick numbers $i_{j}$ so that
\[
{\sf LJm}+(cut)\not\vdash\oplus(\vec{\alp},\oplus(\oplus(\chi^{\fal}(\sig*(0,i_{0})),\bigvee_{j>0}\chi^{\fal}(\sig*(j,i_{j}))),\bet))
.\]
\ref{prp:oplus}.\ref{prp:oplus7}.
Let $\sig$ be a leaf in $T$.
By Proposition \ref{prp:oplus}.\ref{prp:oplus15} we have
${\sf LJm}+(cut)\vdash\chi(T)\lrarw\oplus^{\fal}(\vec{\alp},\chi(\sig))\lrarw\oplus^{\fal}(\vec{\alp},\oplus(\chi(\sig),\bot))$
for some formulas $\vec{\alp}$. 
On the other side the formula $\chi(T^{\prime})$ is obtained from $\chi(T)$ by replacing $\chi(\sig)$ by 
$\oplus(\chi^{\fal}(\sig*(0,i_{0});T),\bigvee_{0<j\leq J_{\supset}}\chi^{\fal}(\sig*(j,i_{j});T)\lor\bigvee_{J_{\supset}<j\leq J_{\fal}}\chi^{\fal}(\sig*(j,i_{j});T))
$
 for the leaf $\sig$ in $T$.
Thus the proposition is seen from Proposition \ref{prp:oplus}.\ref{prp:oplus6}.
\eprf
\\

Supposing the given sequent $S_{0}$ is underivable in ${\sf LJm}+(cut)$, ${\sf LJm}+(cut)\not\vdash S_{0}$,
let us pick
a tree $T\subset TR(S_{0})$
for which the following holds.
Let $T_{n}=\{\sig\in T: lh(\sig)\leq 2n+1\}$.

\benu
\item
 for any $\sig\in T$,
 $g(\sig)\in\{0,\land,\lor\}$,
\item
$\emptyset\in T$, 
and
there exists a unique son $(i)$ of $\emptyset$ in $T$
such that ${\sf LJm}+(cut)\not\vdash S((i))$.
Namely $T_{0}=\{\emptyset,(i)\}$.

Let
\[
\chi(T_{0}):\equiv\chi((i))\equiv(\bigwedge\Gam((i))\supset \bigvee\Del((i)))
\]
Then ${\sf LJm}+(cut)\not\vdash\chi(T_{0})$.

\item
 for any $\sig\in T_{n}$ with $g(\sig)=\lor$,
every son $\sig*(j)\in TR(S_{0})$ is in $T_{n+1}$,
and
there exists a unique son $\sig*(j,i_{j})$ for each $j$.
Namely $T_{n+1}=T_{n}\cup\{\sig*(j),\sig*(j,i_{j}):\sig\in T_{n}, lh(\sig)=2n+1\}$.

Let $lh(\sig)=2n+1$ and assume ${\sf LJm}+(cut)\not\vdash\chi(T_{n})$.
The sons $(j,i_{j})$ are chosen so that 
${\sf LJm}+(cut)\not\vdash\chi(T_{n+1})$.
Such an extension is possible by Proposition \ref{prp:oplus}.\ref{prp:oplus7}.

\eenu

It is clear that $T$ is a subtree of $TR(S_{0})$ such that $g(\sig)\neq 1$ for any $\sig\in T$,
and each $\land$-gate has a unique $\lor$-gate son in $T$.

\blem\label{lem:transfer}

For $\sig\in T$, $\Gam^{\infty}_{\subset_{e}}(\sig;T)\Rarw\Del^{\infty}(\sig;T)$ is 
$FV^{\infty}_{\subset_{e}}(\sig;T)$-analyzed,
and $\Gam^{\infty}_{\subset_{e}}(\sig;T)\cap\Del^{\infty}(\sig;T)$ has no common atomic formula
\elem
\bprf
It is easy to see that for any $\tau\in T$,
$\bot\not\in\Gam^{\times}(\tau)$, and hence $\bot\not\in\Gam^{\infty}_{\subset_{e}}(\sig;T)$.

Since each term over the set $FV^{\infty}_{\subset_{e}}(\sig;T)$ is eventually tested for
universal antecedent formula and existential succedent formula in the extensions
$\tau$ of $\sig$ with $\sig\subset_{e}^{0}\tau$,
$\Gam^{\infty}_{\subset_{e}}(\sig;T)\Rarw\Del^{\infty}(\sig;T)$ is $FV^{\infty}_{\subset_{e}}(\sig)$-saturated.

Suppose that $\alp$ is a common atomic formula in 
$\Gam^{\infty}_{\subset_{e}}(\sig;T)$ and $\Del^{\infty}(\sig;T)$.
Let $\sig_{1}\in T$ be such that 
$g(\sig_{1})\neq 1$, 
$\sig\subset_{e}^{0}\sig_{1}$ and $\alp\in\Del(\sig_{1})^{\times}$.
Also let $\rho,\sig_{0}\in T$ be such that 
$g(\sig_{0})\neq 1$,
 $\rho\subset_{e}\sig$, $\rho\subset_{e}^{0}\sig_{0}$
and $\alp\in\Gam(\sig_{0})^{\times}$.
We see $\rho\subsetneq_{e}\sig$ from Proposition \ref{prp:fanlyzedLJm}, $\alp\not\in\Gam^{\infty}(\sig;T)\cap\Del^{\infty}(\sig;T)$.

We see that $\rho\not\subset_{e}^{0}\sig$,
 otherwise $\sig_{0}$ and $\sig_{1}$ are comparable in the order $\subset_{e}^{0}$,
and one of sequents $S(\sig_{0})$ and $S(\sig_{1})$ is an axiom with $1\in\{g(\sig_{0}),g(\sig_{1})\}\not\ni 1$.
Therefore there exist a $\mu$ and an $i\neq 0$ such that $\rho\subset_{e}\mu$, $\mu*(i)\subset_{e}\sig$ and $g(\mu)=\lor$.
Let $\mu$ be the lowest, i.e., the shortest such sequence.
Then $\rho\subset_{e}^{0}\mu$ and $\mu*(0)\subset_{e}^{0}\sig_{0}$.

We can assume that $\rho=\mu$.
In other words
$\rho$ is the infimum of $\sig_{0}$ and $\sig_{1}$.
\[
\infer[(\mbox{br})]{S(\rho): \Gam\Rarw\Del}
{
\infer*[d_{0}]{\Gam_{0}\Rarw\Del_{0}}
{
S(\sig_{0}): \alp,\Pi_{0}\Rarw\Lam_{0}
}
&
\cdots
&
\infer*[d_{1}]{\Gam_{1}\Rarw\Del_{1}}
{
S(\sig_{1}):\Pi_{1}\Rarw\Lam_{1},\alp
}
&
\cdots
}
\]
Let $2n+1:=\max\{lh(\sig_{0}),lh(\sig_{1})\}$.
Let $a$ be an eigenvariable of an inference rule occurring between $\sig_{1}$ and $\rho$.
If the variable $a$ occurs in the formula $\alp$, then so does in the sequent $S(\sig_{0}): \alp,\Pi_{0}\Rarw\Lam_{0}$,
and this contradicts (\ref{eq:freshvar}).
Hence the variable $a$ does not occur in the formula $\alp$.

There are formulas $\bet$ and $\vphi(\alp)$ 
such that ${\sf LJm}+(cut)\vdash \chi(\rho;T_{n})\lrarw((\alp\land\bet)\supset \vphi(\alp))^{\fal}$,
where $\alp$ occurs in $\vphi(\alp)$ possibly in the scopes of $\lor,\fal$ and in the succedents of $\supset$,
but no universal quantifier binds free variables in $\alp$.
Namely $\vphi(\alp)$ is in the class $\Phi(\alp)$ such that
$\alp\in\Phi(\alp)$, $\vphi(\alp)\in\Phi(\alp)\Rarw \{\bet\lor\vphi(\alp),(\bet\supset\vphi(\alp))^{\fal}\}\subset\Phi(\alp)$
for any $\bet$.
No universal quantifier in $(\bet\supset\vphi(\alp))^{\fal}$ binds free variables in $\alp$
since a universal quantifier in $(\bet\supset\vphi(\alp))^{\fal}$ binds an eigenvariable 
of an inference rule occurring between $\sig_{1}$ and $\rho$.

We see inductively that
${\sf LJm}+(cut)\vdash\alp\supset\vphi(\alp)$ for any $\vphi(\alp)\in\Phi(\alp)$.
Hence ${\sf LJm}+(cut)\vdash \chi(\rho;T_{n})$, and ${\sf LJm}+(cut)\vdash\chi(T_{n})$ by Proposition \ref{prp:oplus}.\ref{prp:oplus1},
${\sf LJm}+(cut)\vdash\alp\lor\bet \supset  \oplus(\alp,\bet)$.
This is a contradiction since $T_{n}$ is chosen so that ${\sf LJm}+(cut)\not\vdash\chi(T_{n})$.
\eprf

\bth\label{th:intcompleteLJm}
Suppose that ${\sf LJm}+(cut)\not\vdash S_{0}$.
Each Kripke model $\langle T,\subset_{e},D_{T},I_{T}\rangle$ falsified the given sequent $S_{0}$,
no matter which tree $T$ is chosen.

Hence ${\sf LJm}+(cut)$ is intuitionistically complete in the sense that
any intuitionistically valid sequent is derivable in ${\sf LJm}+(cut)$.
\end{theorem}
\bprf
$T$ enjoys the three conditions in Proposition \ref{prp:countermodel}.
The third condition follows from Lemma \ref{lem:transfer}.
Hence for $S_{0}=S(\emptyset)=(\Gam(\emptyset)\Rarw\Del(\emptyset))$,
$\emptyset\models\bigwedge\Gam(\emptyset)$ and $\emptyset\not\models\bigvee\Del(\emptyset)$
in the Kripke model 
 $\langle T,\subset_{e},D_{T},I_{T}\rangle$ defined from the tree $T$.
 \eprf

\subsection{Transfer}\label{subsec:transfer}
It may be the case that for a $\sig\in TR(S_{0})$,
$\Gam^{\infty}_{\subset_{e}}(\sig;TR(S_{0}))\cap\Del^{\infty}(\sig;TR(S_{0}))$ has a common atomic formula,
and we need to transform $TR(S_{0})$.

We consider a transformation of deductions inspired by the transfer rule in \cite{Kripke}.
The tree $TR(S_{0})$ of deductions is transformed to another tree $TTR(S_{0})$ in which there is no $\sig$
such that $\Gam^{\infty}_{\subset_{e}}(\sig;T)\cap\Del^{\infty}(\sig;T)$ has a common atomic formula,
where $T$ is a subtree of $TTR(S_{0})$ such that each $\land$-gate has a unique $\lor$-gate son in $T$.

Let us modify the inference rule $({\rm br})$ by combining a weakening as follows.

{\tiny
\[
\infer[(\mbox{wbr})]{\Pi,\Gam,\Gam_{\supset},\Gam_{\fal}\Rarw\Lam,\Del_{\exi},\Del,\Del_{\supset}^{\circ},\Del_{\fal}^{\circ}}
{
\Pi,\Gam,\Gam_{\supset},\Gam_{\fal}^{\circ}\Rarw \Lam,\Del_{\exi}^{\circ},\Del,\Del_{\supset}^{\circ},\Del_{\fal}^{\circ}
&
\{\gam^{\circ},\Gam,\Gam_{\supset}^{\circ},\Gam_{\fal}^{\circ}\Rarw\del^{\circ}: (\gam\supset\del)^{\circ}\in\Del_{\supset}^{\circ}\}
& 
\{\Gam,\Gam_{\supset}^{\circ},\Gam_{\fal}^{\circ}\Rarw\gam(a)^{\circ}:(\fal y\,\gam(y))\in\Del_{\fal}^{\circ}\}
}
\]
}
where $\Pi$ and $\Lam$ are arbitrary cedents.

Moreover let us introduce the following inference rule.
\[
\infer[(\circ)]{\Gam,\Gam_{\fal}\Rarw\Del_{\exi},\Del}
{\Gam,\Gam_{\fal},\Gam_{\fal}^{\circ}\Rarw\Del_{\exi}^{\circ},\Del_{\exi},\Del}
\]
where $\Gam_{\fal}$ is the set of unmarked universal formulas in the antecedent of lower sequent,
$\Del_{\exi}$ a set of unmarked existential formulas in the succedent.

Let $T$ be a (finite or infinite) labelled tree of deductions in which the inference rules 
$(\mbox{wbr}), 
(\circ)$ may occur besides
inference rules in {\sf LJm}, and $S(\sig)$ is a sequent for $\sig\in T$.
Suppose that the tree $T$ of deductions enjoys the condition (\ref{eq:freshvar})
on the eigenvariables.

For $\sig_{0},\sig_{1}\in T$, we say that a pair $(\sig_{0},\sig_{1})$ is a \textit{transferable pair}
if the following conditions are met, cf. the proof of Lemma \ref{lem:transfer}:
\benu
\item 
The antecedent of $S(\sig_{0})$ and the succedent of $S(\sig_{1})$ have a common marked atomic formula 
$\alp^{\circ}\equiv (R(t_{1},\ldots,t_{n}))^{\circ}$.
This means that $S(\sig_{0})$ is a sequent $\alp^{\circ},\Pi_{0}\Rarw\Lam_{0}$,
and $S(\sig_{1})$ is $\Pi_{1}\Rarw\Lam_{1},\alp^{\circ}$ for some cedents $\Pi_{i},\Lam_{i}$.

\item
Each inference rule at $\sig_{0}$
 and at the infimum $\rho$ of $\sig_{0}, \sig_{1}$ is a $({\rm wbr})$.

\item
There is no non-invertible upper sequent from $S(\sig_{0})$ to the sequent $S(\rho*(0))$
in the deduction, i.e.,
$\rho\subset_{e}^{0}\sig_{0}=\rho*\kap_{0}$ for some $\kap_{0}=(0)*\kap_{0}^{\prime}$.

And $\rho\subsetneq_{e}\sig_{1}=\rho*\kap_{1}$ with $\rho\not\subset_{e}^{0}\sig_{1}$
for some $\kap_{1}=(i)*\kap_{1}^{\prime}$ and $i\neq 0$.

\eenu

All of transferrable pairs have to be removed from the constructed deduction in proof search
to avoid an inconsistency in the definition of Kripke models,
cf. Proposition \ref{prp:countermodel}
and Lemma \ref{lem:transferred}.

For this, let us transform the tree $TR(S_{0})$ of deductions to another tree $TTR(S_{0})$.
The gate $g(\sig)$ is changed simultaneously in the transformations.

\bdf\label{df:gatechange}
{\rm
Let $g(\tau,0)=g(\tau)$ for $\tau\in TR(S_{0})$.
A gate $g(\tau,n)$ is computed as follows.
Consider the case when $g(\tau,n)=\lor$.
If there exists an $i$ such that $g(\tau*(i),n+1)=1$,
then $g(\tau,n+1)=1$.
If $g(\tau*(i),n+1)=0$ for any $i$, then $g(\tau,n+1)=0$.
Second consider the case when $g(\tau,n)=\land$.
If there exists an $i$ such that $g(\tau*(i),n+1)=0$,
then $g(\tau,n+1)=0$.
If $g(\tau*(i),n+1)=1$ for any $i$, then $g(\tau,n+1)=1$.
In all other cases the gate is unchanged, $g(\tau,n+1)=g(\tau,n)$.
}
\edf

Let $TTR(S_{0},0)=TR(S_{0})$ with the same labeling functions $S(\sig,0)=S(\sig)$,
$d(\sig,0)=d(\sig)$ and $g(\sig,0)=g(\sig)$.
Suppose that $TTR(S_{0},n)$ has been defined.
Pick a minimal transferable pair $(\sig_{0},\sig_{1})$ in $TTR(S_{0},n)$, 
where the minimality means that the length $lh(\rho)$ of the infimum
$\rho$ of $\sig_{0}$ and $\sig_{1}$ is minimal.
Let $\sig_{i}=\rho*\kap_{i}$ for $i=0,1$.

\[
\infer[(\mbox{wbr})]{\rho:\Pi,\Gam,\Gam_{\supset},\Gam_{\fal}\Rarw \Del_{\exi},\Del,\Del_{\supset}^{\circ},\Del_{\fal}^{\circ},\Lam}
{
\infer*[d_{0}]{\Pi,\Gam,\Gam_{\supset},\Gam_{\fal}^{\circ}\Rarw \Del_{\exi}^{\circ},\Del,\Del_{\supset}^{\circ},\Del_{\fal}^{\circ},\Lam}
{
 \infer*{\rho*\kap_{0}: \alp^{\circ},\Pi_{0}\Rarw\Lam_{0}}{}
}
&
\cdots
&
\infer*[d_{1}]{\Gam_{1},\Gam_{\supset}^{\circ},\Gam_{\fal}^{\circ}\Rarw\del_{1}}
{
\infer*{\rho*\kap_{1}:\Pi_{1}\Rarw\Lam_{1},\alp^{\circ}}{}
}
}
\]
where $d_{0}$ denotes a deduction of the leftmost continued sequent up to $\rho*\kap_{0}$,
and $d_{1}$ a deduction of a non-invertible upper sequent $\Gam,\Gam_{\supset}^{\circ},\Gam_{\fal}^{\circ}\Rarw\del_{1}$ up to $\rho*\kap_{1}$.

Then let us combine these two deductions, and 
 transfer the current tree $TTR(S_{0},n)$ by the pair $(\sig_{0},\sig_{1})$ to get $TTR(S_{0},n+1)$.
\[
\infer[(\circ)]{\rho:\Pi,\Gam,\Gam_{\supset},\Gam_{\fal}\Rarw\Del_{\exi},\Del,\Del_{\supset}^{\circ},\Del_{\fal}^{\circ},\Lam}
{
\infer*[\Gam_{2}*d_{0}*\Del_{2}]{\Pi,\Gam,\Gam_{\supset},\Gam_{\fal}^{\circ},\Gam,\Gam_{\supset},\Gam_{\fal} \Rarw\Del_{\exi},\Del_{\exi}^{\circ},\Del,\Del_{\supset}^{\circ},\Del_{\fal}^{\circ},\Lam}
{
  \infer[(\mbox{wbr})]{\rho*\kap_{0}:\alp^{\circ},\Pi_{0},\Pi,\Gam,\Gam_{\supset},\Gam_{\fal}\Rarw\Del_{\exi},\Del,\Del_{\supset}^{\circ},\Del_{\fal}^{\circ},\Lam,\Lam_{0}}
  {
  \infer*[ \alp^{\circ}\Pi_{0}*d_{0}^{\prime}*\Lam_{0}]{ \alp^{\circ},\Pi_{0},\Pi,\Gam,\Gam_{\supset},\Gam_{\fal}^{\circ}\Rarw\Del_{\exi}^{\circ},\Del,\Del_{\supset}^{\circ},\Del_{\fal}^{\circ},\Lam,\Lam_{0}}
   {
     \infer*{\rho*\kap_{0}*\kap_{0}: \alp^{\circ},\Pi_{0},\alp^{\prime\circ},\Pi^{\prime}_{0}\Rarw\Lam^{\prime}_{0},\Lam_{0}}{}
    }
 &
 \cdots
 &
  \infer*[\alp^{\circ}*d_{1}]{\alp^{\circ},\Gam_{1},\Gam_{\supset}^{\circ},\Gam_{\fal}^{\circ}\Rarw\del_{1}}
  {
   \rho*\kap_{0}*\kap_{1}:\alp^{\circ},\Pi_{1}\Rarw\Lam_{1},\alp^{\circ}
   }
   }
 }
}
\]
where $\rho*\kap_{0}*\kap_{1}:\alp^{\circ},\Pi_{1}\Rarw\Lam_{1},\alp^{\circ}$ is an axiom, 
$\alp^{\circ}*d_{1}$ is a deduction obtained from $d_{1}$ by appending $\alp^{\circ}$
to antecedents.
$\Gam_{2}*d_{0}*\Del_{2}$ is a deduction obtained from $d_{0}$ by
appending $\Gam_{2}=\Pi\cup\Gam\cup\Gam_{\supset}\cup\Gam_{\fal}$ to antecedents and 
$\Del_{2}=\Del_{\exi}\cup\Del\cup\Del_{\supset}^{\circ}\cup\Del_{\fal}^{\circ}\cup\Lam$ to succedents.

$\alp^{\circ}\Pi_{0}*d_{0}^{\prime}*\Lam_{0}$ is a deduction obtained from $d_{0}$
by appending $\alp^{\circ}\cup\Pi_{0}$ to the antecedents, $\Lam_{0}$ to the succedents,
 and renaming eigenvariables from ones 
in $d_{0}$ for the condition (\ref{eq:freshvar}) on eigenvariables.
$\alp^{\prime\circ}$ ($\Pi_{0}^{\prime}$) [$\Lam_{0}^{\prime}$] is obtained from 
$\alp^{\circ}$ ($\Pi_{0}$) [$\Lam_{0}$]
by renaming free variables which are introduced as eigenvariables of $(\exi\Rarw)$
in $d_{0}$, resp.

The above deduction is said to be \textit{transferred} by the pair $(\sig_{0},\sig_{1})$
(and the common atomic formula $\alp^{\circ}$).

Let $g(\sig,n+1)\in\{0,1,\lor,\land\}$ be the gate labeling function for the transferred deduction.
In the transferred deduction, each sequent up to 
$\rho*\kap_{0}:\alp^{\circ},\Gam,\Gam_{\supset},\Gam_{\fal}\Rarw\Del_{\exi},\Del,\Del_{\supset}^{\circ},\Del_{\fal}^{\circ}$ receives the same label $\mu$ 
and the same gate $g(\mu,n+1)=g(\mu,n)$ 
as in the original deduction, e.g., $g(\rho,n+1)=\lor=g(\rho,n)$ and $g(\rho*(0),n+1)=\land=g(\rho*(0),n)$.
On the other side in the above part of 
$\rho*\kap_{0}:\alp^{\circ},\Gam,\Gam_{\supset},\Gam_{\fal}\Rarw\Del_{\exi},\Del,\Del_{\supset}^{\circ},\Del_{\fal}^{\circ}$,
insert $\kap_{0}$ to each label where $\sig_{0}=\rho*\kap_{0}$.
$g(\rho*\kap_{0}*\kap_{0},n+1)=g(\rho*\kap_{0},n+1)=\lor=g(\rho*\kap_{0},n)$.
In the transferred deduction $TTR(S_{0},n+1)$ the gate of the node $\rho*\kap_{0}*\kap_{1}$ becomes $1$,
$g(\rho*\kap_{0}*\kap_{1},n+1)=1$, 
since $S(\rho*\kap_{0}*\kap_{1},n+1)$ is an axiom.
Below the axiom $\rho*\kap_{0}*\kap_{1}$, modify the values of gates as in Definition \ref{df:gatechange}.
Some gates might receive the value $1$ by this modifications.

Let us check the conditions on multi-succedents and eigenvariables.
Since there is no non-invertible upper sequent from $\rho$ to $\rho*\kap_{0}$, we can append formulas in $\Lam_{0}$ and
formulas in
$\Del_{2}=\Del_{\exi}\cup\Del\cup\Del_{\supset}^{\circ}\cup\Del_{\fal}^{\circ}$ to the succedents.
Next consider a free variable $a$ occurring in $\alp$.
Assume that $a$
is introduced as an eigenvariable of an $(L\exi)$ or a $(\mbox{wbr})$(, 
which subsumes $(R\fal)$)
between $\rho$ and $\rho*\kap_{1}$ in $d_{1}$.
But the variable $a$ occurs in $\rho*\kap_{0}$, and this is not the case by (\ref{eq:freshvar})
since $\rho*\kap_{0}$ is not above the inference rule nor to the right of it.
Furthermore free variables occurring in $\Gam_{2}\cup\Del_{2}$ is distinct from
 eigenvariables of $(L\exi)$ occurring in $d_{0}$
since $\rho$ is below $d_{0}$, (\ref{eq:freshvar}).

Iterate the transformations to yield a tree of deductions $TTR(S_{0})$, in which there is no transferable pair.

This ends the construction of $TTR(S_{0})\subset{}^{<\ome}\ome$,
where 
\[
TTR(S_{0})=\liminf_{n\to\infty}TTR(S_{0},n)=\bigcup_{n}\bigcap_{m\geq n}TTR(S_{0},m).
\]
$S(\sig,n),d(\sig,n),g(\sig,n)$ are labeling functions of sequents, deductions and gates for
the nodes $\sig\in TTR(S_{0},n)$.

Let
for $\sig\in TTR(S_{0})$
\beqnarrs
\bar{S}(\sig) & = & \lim_{n\to\infty}S(\sig,n)
\\
\bar{d}(\sig) & = & \lim_{n\to\infty}d(\sig,n)
\\
\bar{g}(\sig) & = & \lim_{n\to\infty}g(\sig,n)
\eeqnarrs

These limits exist since for any $k$ 
there exists an $n_{k}$ such that for any $m$ and any transferable pair 
$(\sig_{0},\sig_{1})$ in $TTR(S_{0},n_{k}+m)$,  $lh(\rho)>k$ holds for the infimum $\rho$ of $\sig_{0},\sig_{1}$.
Then $S(\sig,n_{k}+m)=S(\sig,n_{k})$ and $d(\sig,n_{k}+m)=d(\sig,n_{k})$ for any $m$ and any $\sig$ with $lh(\sig)\leq k$.
Furthermore for such a $\sig$, the series
$\{g(\sig,n_{k}+m)\}_{m}$ changes the values monotonically, $g(\sig,n_{k}+m)\leq g(\sig,n_{k}+m+l)$
 in the order $\lor,\land<0<1$ on $\{0,1,\lor,\land\}$.
Hence the limit $\bar{g}(\sig)$ exists, too.

\blem\label{lem:transferred}
Let $\emptyset\neq T\subset TTR(S_{0})$ be a subtree such that $\fal\sig\in T[\bar{g}(\sig)\neq 1]$ and
each $\land$-gate has a unique $\lor$-gate son in $T$, i.e., 
$\fal\sig\in T\exi! i[\bar{g}(\sig)=\land \Rarw \bar{g}(\sig*(i))=\lor \spand \sig*(i)\in T]$.
Then
for $\sig\in T$, 
there is no common atomic formula in
$\Gam^{\infty}_{\subset_{e}}(\sig;T)\cap\Del^{\infty}(\sig;T)$.
\elem
\bprf
Suppose that $\alp$ is a common atomic formula in 
$\Gam^{\infty}_{\subset_{e}}(\sig;T)$ and 
$\Del^{\infty}(\sig;T)$.
Let 
 $\sig_{1}\in T$ be such that  $\alp\in\Del(\sig_{1})^{\times}$ and
 $\sig\subset_{e}^{0}\sig_{1}$.

First consider the case when $\alp\in\Gam(\tau_{0})^{\times}$ for some $\tau_{0}\subset_{e}\sig$.
Then $\alp\in\Gam(\sig_{1})^{\times}$.
 This means that $S(\sig_{1})$ is an axiom, and $\bar{g}(\sig_{1})=1$.

Let $\rho,\sig_{0}\in T$ be such that $\bar{g}(\sig_{0},lh(\sig_{0}))=\lor$, 
$\rho\subset_{e}\sig$, $\rho\subset_{e}^{0}\sig_{0}$
and $\alp\in\Gam(\sig_{0})^{\times}$.
We see $\rho\subsetneq_{e}\sig$ from Proposition \ref{prp:fanlyzedLJm}, 
$\alp\not\in\Gam^{\infty}(\sig;T)\cap\Del^{\infty}(\sig;T)$.
Furthermore $\rho\not\subset_{e}^{0}\sig$ otherwise $\sig_{0}$ and $\sig_{1}$ are comparable in the order $\subset_{e}^{0}$,
and one of sequents $S(\sig_{0})$ and $S(\sig_{1})$ is an axiom with 
$1\in\{\bar{g}(\sig_{0},lh(\sig_{0})),\bar{g}(\sig_{1},lh(\sig_{1}))\}\not\ni 1$.
Therefore there exist a $\tau$ and an $i\neq 0$ such that $\rho\subset_{e}\tau$, $\tau*(i)\subset_{e}\sig$ and $\bar{g}(\tau,lh(\tau))=\lor$.
Let $\tau$ be the lowest, i.e., the shortest such sequence.
Then $\rho\subset_{e}^{0}\tau$ and $\tau*(0)\subset_{e}^{0}\sig_{0}$.
This means that $(\sig_{0},\sig_{1})$ is a transferable pair.
Such a pair has been removed from $TTR(S_{0},n)$ for an $n$, and from $TTR(S_{0})$.
Hence this is not the case.
\eprf

\bth\label{th:schuttedichotomy}{\rm (Sch\"utte's dichotomy)}\\
In $TTR(S_{0})$, $\bar{g}(\emptyset)=1$ iff ${\sf LJm}\vdash S_{0}$.
\end{theorem}
\bprf
If $\bar{g}(\emptyset)=1$, then it is plain to see that ${\sf LJm}\vdash S_{0}$.

In what follows assume $\bar{g}(\emptyset)\neq 1$.
Then $\bar{g}(\emptyset)\in\{0,\land\}$.
Extract a subtree $T\subset TTR(S_{0})$ as follows.
First $\emptyset\in T$.
The nodes $\sig\in T$ with $\bar{g}(\sig)=0$ are leaves in $T$.
Suppose $\sig\in T$ has been chosen so that $\bar{g}(\sig)=\land$.
Then in the deduction $d(\sig)$, pick a leaf $\sig*(j)$ such that $\bar{g}(\sig*(j))\neq 1$.
If $\bar{g}(\sig*(j))=0$, then we would have $\bar{g}(\sig)=0$.
Hence $\bar{g}(\sig*(j))=\lor$, and $\bar{g}(\sig*(j,i))\neq 1$ for any son $\sig*(j,i)$.
Moreover there exists an $i$ such that $\bar{g}(\sig*(j,i))\neq 0$.
Otherwise we would have $\bar{g}(\sig*(j))=0$.
Let $\sig*(i),\sig*(j,i)\in T$ for any $i$ such that $\bar{g}(\sig*(j,i))=\land$.

Then $T$ enjoys the three conditions in Proposition \ref{prp:countermodel}.
The third condition follows from Lemma \ref{lem:transferred}.
Hence for $S_{0}=S(\emptyset)=(\Gam(\emptyset)\Rarw\Del(\emptyset))$,
$\emptyset\models\bigwedge\Gam(\emptyset)$ and $\emptyset\not\models\bigvee\Del(\emptyset)$
in the Kripke model 
 $\langle T,\subset_{e},D_{T},I_{T}\rangle$ defined from the tree $T$.
 Therefore ${\sf LJm}\not\vdash S_{0}$.
 \eprf

\bcor
Cut inference
\[
\infer[(cut)]{\Gam\Rarw\Del}
{\Gam\Rarw \alp,\Del
&
\Gam,\alp\Rarw\Del
}
\]
is permissible in {\sf LJm}.
\ecor


\begin{thebibliography}{99}

\bibitem{BussIemhoff}
S. Buss and R. Iemhoff, 
The depth of intuitionistic cut free proofs, 
unpublished, available online, 2003.


\bibitem{Dowek} G. Dowek and J. Ying, 
Eigenvariables, bracketing and the decidability of positive minimal predicate logic,
Theoret. Comput. Sci., 360(2006), 193-208.

\bibitem{Dyckhoff}R. Dyckhoff,
Intuitionistic decision procedures since Gentzen,
in: Advances in Proof Theory,
R. Kahle, T. Strahm and T. Studer, eds.,
Birkh\"auser, 2016, pp. 245-267.


\bibitem{Hudelmaier}
J. Hudelmaier,
An $O(n\log n)$-space decision procedure for intuitionistic propositional logic,
J. Logic Computat. 3(1993), 63-75.


\bibitem{Kripke}S. Kripke, 
Semantic analysis of intuitionistic logic, 
in: J. Crossley and M. A. E. Dummett, eds., Formal Systems and Recursive Functions, 
North-Holland, Amsterdam (1965), pp. 92-130.

\bibitem{Ladner}
R. E. Ladner,
The computational complexity of provability in systems of modal propositional logic,
Siam. J. Comput., 6(1977), 467-480.

\bibitem{Maehara}S. Maehara,
Eine Darstellung der intuitionistischen Logik in der klassischen,
Nagoya Math. J. 7(1954), 45-64.

\bibitem{MintsSol}
G. Mints, 
Solvability of the problem of deducibility in LJ for a class of formulas not containing negative occurrences of quantifiers, 
Proc. Steklov Inst. Math. 98(1968), 135-145.


\bibitem{Mintsmodal}G. Mints,
A Short Introduction to Modal Logic, CSLI Lect. Notes 30, Stanford,1992.

\bibitem{Mintsint}G. Mints, 
A Short Introduction to Intuitionistic Logic, Kluwer Academic/Plenum Publishers, New York, 2000.


\bibitem{Mintsdraft}G. Mints, 
Structure of tableau derivations and cut-elimination for intuitionistic logic, draft.


\bibitem{PintoDyckhoff}
L. Pinto and R. Dyckhoff,
Loop-free construction of counter models for intuitionistic propositional logic,
in Symposia Gaussiana, Conf. A, Eds.: Behara/Fritsch/Lintz, Walter de Gruyter (1995), pp. 225-232.

\bibitem{Schuette56}K. Sch\"utte, Ein System des verkn\"upfenden Schliessens,
Archiv math. Logik Grundlagenforsch. (1956), 55-67.


\bibitem{Schuette60}K. Sch\"utte, Syntactical and semantical properties of simple type theory,
Jour. Symb. Logic 25(1960), 305-326.

\bibitem{TroSch}A. S. Troelstra and H. Schwichteberg,
Basic Proof Theory, second edition, Cambridge UP, 2000.

\end{thebibliography}
\end{document}